\input ./style/arxiv-general.cfg
\documentclass[MSNbibl,number,citesort,seceqn,secthm,dvips]{arxbj}
\makeatletter
   \@ifpackageloaded{graphicx}{}{\usepackage{graphicx}}
\makeatother

\volume{23}
\issue{1}
\pubyear{2017}
\firstpage{159}
\lastpage{190}
\doi{10.3150/15-BEJ740}
\docsubty{FLA}

\makeatletter

\newcommand{\eqref}[1]{(\ref{#1})}
\newtheorem{coro}{Corollary}[section]
\newtheorem{lem}{Lemma}[section]

\newremark{rem}{Remark}
\newproclaim{defi}{Definition}[section]
\newproclaim{exm}{Example}[section]

\makeatother

\begin{document}
\begin{frontmatter}

\title{Mixed domain asymptotics for a stochastic process model with
time trend and measurement error}
\runtitle{Mixed domain asymptotics}

\begin{aug}
\author[A]{\inits{C.-H.}\fnms{Chih-Hao}~\snm{Chang}\thanksref{A}\ead[label=e1]{jhow@nuk.edu.tw}},
\author[B]{\inits{H.-C.}\fnms{Hsin-Cheng}~\snm{Huang}\thanksref{B,e2}\ead[label=e2,mark]{hchuang@stat.sinica.edu.tw}}
\and
\author[B,C]{\inits{C.-K.}\fnms{Ching-Kang}~\snm{Ing}\thanksref{B,C,e3}\corref{}\ead[label=e3,mark]{cking@stat.sinica.edu.tw}\ead[label=e4]{cking1@ntu.edu.tw}}
\address[A]{Institute of Statistics, National University of Kaohsiung, Kaohsiung 811, Taiwan.\\ \printead{e1}}
\address[B]{Institute of Statistical Science, Academia Sinica, Taipei 115, Taiwan.\\ \printead{e2,e3}}
\address[C]{Department of Economics, National Taiwan University, Taipei 106, Taiwan.\\ \printead{e4}}
\end{aug}

%
\received{\smonth{5} \syear{2014}}
%
\revised{\smonth{5} \syear{2015}}

%
\begin{abstract}
We consider a stochastic process
model with time trend and measurement error.
We establish consistency and derive the limiting
distributions of the maximum likelihood (ML) estimators of the
covariance function parameters under a general asymptotic framework,
including both the fixed domain and the increasing domain
frameworks, even when the time trend model
is misspecified or its complexity increases with the
sample size. In particular,
the convergence rates of the ML estimators are thoroughly characterized
in terms of
the growing rate of the domain and the degree of model
misspecification/complexity.
\end{abstract}

%
\end{frontmatter}

\section{Introduction}\label{sec1}

Learning the covariance structure of a stochastic
process from data is a fundamental prerequisite for problems such as
prediction, classification and control.
For example, to do prediction for an Ornstein--Uhlenbeck (OU) (Uhlenbeck
and Ornstein \cite{Uhlenbeck})
process $\eta(s)$, $s \in[0,1]$ with mean 0
and covariance function
%
\begin{eqnarray}\label{expcovfun1dim}
&&\operatorname{cov}\bigl(\eta(s_{1}),\eta(s_{2})\bigr) =
\sigma_{0,\eta}^2\exp\bigl(-\kappa_{0}|s_{1}-s_{2}|\bigr),
\end{eqnarray}
where $\sigma_{0,\eta}^2, \kappa_{0}>0$ are unknown, Ying \cite{Ying1991}
proposed the maximum likelihood (ML) estimators for
$\sigma_{0,\eta}^2$ and $\kappa_{0}$ based on discrete observations
$\eta(s_{1}), \ldots, \eta(s_{n})$ with $0\leq
s_{1}<\cdots<s_{n}\leq1$, and established the root-$n$ consistency
of the corresponding ML estimator for $\sigma_{0,\eta}^2\kappa_{0}$.
Note that since the probability measures induced by two OU processes
are absolutely continuous with respect to each other if and only if
their $\sigma_{0,\eta}^2\kappa_{0}$ values are equal (Ibragimov and
Rozanov \cite{Ibragimov}), the parameters in (\ref{expcovfun1dim}) are
asymptotically identifiable up to
$\sigma_{0,\eta}^2\kappa_{0}$.
However, when the OU process is subject to measurement error,
the so-called ``nugget'' effect (see, for example, Cressie \cite{Cressie})
may deteriorate the performance of the ML estimators.
In particular, Chen, Simpson and Ying \cite{Chen} showed that
the ML estimator for $\sigma_{0,\eta}^2\kappa_{0}$ becomes
fourth-root-$n$ consistent, depicting
the effect of measurement error in estimating
the exponential covariance parameters
in (\ref{expcovfun1dim}).
On the other hand, they also proved that
the ML estimator of the measurement-error variance
has the usual root-$n$ consistency.
In fact, a similar phenomenon can also be found in a driftless Brownian
motion (BM) process
with measurement error. Let this error-contaminated process be denoted
by $y(t), t \in[0,1]$.
Having observed $y(0), y(1/n), \ldots, y(1)$,
Stein \cite{Stein1990} showed that
a modified ML (MML) estimator of the ratio of the variance of the increments
of the BM process to that of measurement error is only
fourth-root-$n$ consistent,
whereas the corresponding MML estimator
of
the measurement-error variance
still remains root-$n$ consistent.
Similar asymptotic results for the
ML estimators of the two variances
have also been established by
A\"{\i}t-Sahalia, Mykland and Zhang \cite{AMZ}.

In this article, we shall superimpose a time trend (regression) term on
an OU process
with measurement error in order to accommodate a broader range of applications.
Specifically, we propose the following model for a real-valued
stochastic process
$\{Z(s); s \in D \subset\mathbb{R}\}$:
%
\begin{eqnarray}\label{setupfullmodel}
&& Z(s)=\beta_0+\sum_{j=1}^{p}
\beta_jx_j(s)+\eta(s)+\epsilon(s),
\end{eqnarray}
where ${\mathbf x}(s)=(x_{1}(s), \ldots, x_{p}(s))^{\prime}$
is a $p$-dimensional time trend vector, $\eta(s)$ is a zero-mean OU process
with covariance function defined in (\ref{expcovfun1dim}),
$\epsilon(s)$ is a zero-mean Gaussian measurement error
with $\mathrm{E}(\epsilon(s)\epsilon(t))=\theta_{0, 1}I_{\{s=t\}}$ for
some unknown $\theta_{0, 1}>0$,
$\bolds\beta=(\beta_{0}, \beta_{1}, \ldots, \beta_{p})^{\prime}$
is a ${(p+1)}$-dimensional constant vector,
and $\{{\mathbf x}(s)\}$, $\{\eta(s)\}$ and $\{\epsilon(s)\}$ are independent.
In a computer experiment,
$\eta(s)$ in \eqref{setupfullmodel}
can be used to describe
the systematic departure of the response $Z(s)$
from the linear model
$\beta_0+\sum_{j=1}^{p}\beta_jx_j(s)$
and $\epsilon(s)$
denotes\vspace*{1pt} the measurement error.
For more details, we refer the reader to
Sacks, Schiller and Welch \cite{SSW}
and Ying \cite{Ying1991}.
Model \eqref{setupfullmodel}
can also be applied to
one-dimensional geostatistical modeling
and $\eta(\cdot)$ therein corresponds to a commonly
used exponential covariance model;
see Ripley \cite{Rip} and Cressie \cite{Cressie} for numerous examples.
Denote the true time trend
by
%
\begin{eqnarray}\label{setupfullmodelF}
&& \mu_{0}(s)=Z(s)-\eta(s)-\epsilon(s),
\end{eqnarray}
where $\{\mu_{0}(s)\}$ is independent of
$\{\eta(s)\}$ and $\{\epsilon(s)\}$,
and define ${\mathbf x}_{0}(s)=(1, {\mathbf x}^{\prime}(s))^{\prime}$.
The time trend model $\bolds\beta^{\prime}{\mathbf x}_{0}(s)$ in
(\ref{setupfullmodel}) is said to be correctly specified if
%
\begin{eqnarray}\label{correct}
&& \mu_{0}(s)=\bolds\beta^{\prime}{\mathbf x}_{0}(s)\qquad
\mbox{for some } \bolds\beta\in\mathbb{R}^{p+1},
\end{eqnarray}
and misspecified otherwise.
In this article, we shall allow $\bolds\beta^{\prime}{\mathbf
x}_{0}(s)$ to be misspecified,
which further increases the flexibility of model \eqref{setupfullmodel}.
However, a misspecified time trend will usually create extra challenges
in estimating
covariance parameters.
This motivates us to ask how
the ML estimators of the covariance parameters
in model \eqref{setupfullmodel} perform
when the corresponding time trend model is subject to misspecification.

To facilitate exposition, we assume in the sequel that $D=[0,n^\delta]$
for some
$\delta\in[0,1)$, and the data are observed regularly at
$s_i=in^{-(1-\delta)}$, $i=1,\ldots,n$.
In addition, we also allow that the number of regressors (model complexity)
$p=p_{n}$ grows to infinity in order to reduce the model bias.
When $\delta=0$, the domain $D=[0,1]$ has been considered
by the aforementioned authors, and the setup is called fixed domain asymptotics.
On the other hand, when $\delta>0$, the domain $D$ grows
to infinity as $n\rightarrow\infty$ with a faster growing rate for a larger
$\delta$ value, and the setup is referred to as the increasing domain
asymptotics, even though the minimum inter-data distance
$n^{-(1-\delta)}$ goes to zero.
This is different from the increasing domain setup considered by Zhang
and Zimmerman \cite{Zhang2005}, in which the minimum distance between
sampling points is bounded away from zero.
By incorporating both fixed and increasing domains,
our {\it mixed} domain asymptotic framework
enables us to explore the interplay between
the model misspecification/complexity and the growing rate of $D$ on
the asymptotic behaviors of the ML estimators,
thereby leading to an intriguing answer to the above question.

Re-parameterizing (\ref{expcovfun1dim})
by $\theta_{0, 2}=\sigma_{0,\eta}^2\kappa_{0}$ and $\theta_{0,3}=\kappa_{0}$,
the covariance parameter vector in model (\ref{setupfullmodel}) can be
written as
$\bolds{\theta}_{0}=(\theta_{0,1},\theta_{0,2},\theta_{0,3})'$.
Let $\Theta$, the parameter space, be a compact set in
$(0,\infty)^3$ and suppose
$\bolds{\theta}_{0} \in\Theta$.
Based on model (\ref{setupfullmodel})
and observations $({\mathbf x}^{\prime}(s_i), Z(s_i)), i=1, \ldots,  n$,
we estimate $\bolds{\theta}_{0}$
using the ML estimator
$\hat{\bolds\theta}$, which satisfies
\begin{eqnarray*}
&& \ell(\hat{\bolds\theta})= \sup_{\bolds\theta=(\theta_{1}, \theta_{2},
\theta_{3})^{\prime} \in\Theta} \ell(\bolds\theta),
\end{eqnarray*}
where
%
\begin{eqnarray}
\ell(\bolds\theta) &= & -\tfrac{1}{2}n\log(2\pi)-\tfrac{1}{2}\log\det
\bigl({\bolds\Sigma}({\bolds\theta})\bigr)
\nonumber
\\[-8pt]
\label{loglikefuntrue}
\\[-8pt]
\nonumber
&&{}-\tfrac{1}{2}\mathbf{Z}'\bigl(\mathbf{I}-\mathbf{M}(\bolds
\theta)\bigr)'\bolds\Sigma^{-1} ({\bolds\theta}) \bigl(
\mathbf{I}-\mathbf{M}(\bolds\theta)\bigr)\mathbf{Z},
\end{eqnarray}
is known as the profile log-likelihood function, in which
${\mathbf Z}=(Z(s_1),\ldots,Z(s_n))'$,
%
\begin{eqnarray}\label{Sigma}
&& \bolds\Sigma(\bolds\theta)=\bolds\Sigma_\eta(\bolds\theta)+
\theta_1\mathbf{I},
\end{eqnarray}
with
\begin{eqnarray*}
&& \bolds\Sigma_\eta(\bolds\theta)= \biggl(\frac{\theta_{2}}{\theta
_{3}}\exp\bigl(-
\theta_{3}|s_{i}-s_{j}|\bigr) \biggr)_{1\leq i,j\leq n},
\end{eqnarray*}
and
%
\begin{eqnarray}\label{projectionmatrix}
&& \mathbf{M}(\bolds\theta) = \mathbf{X}\bigl({\mathbf X}'{\bolds
\Sigma}^{-1}(\bolds\theta){\mathbf X}\bigr)^{-1} {\mathbf
X}'{\bolds\Sigma}^{-1}(\bolds\theta),
\end{eqnarray}
with ${\mathbf X}=(\mathbf{x}_{0}(s_1), \ldots, \mathbf
{x}_{0}(s_n))^{\prime}$
being full rank almost surely (a.s.).
It is not difficult to show that the ML estimator of $\bolds{\beta}$ is given
by $\hat{\bolds{\beta}}(\hat{\bolds{\theta}})$, where
\[
\hat{\bolds\beta}(\bolds{\theta})=\bigl(\mathbf{X}'\bolds
\Sigma^{-1}(\bolds\theta)\mathbf{X}\bigr)^{-1}
\mathbf{X}'{\bolds\Sigma}^{-1}(\bolds\theta)\mathbf{Z}.
\]
However, since model (\ref{setupfullmodel}) can be misspecified,
investigating the asymptotic properties of $\hat{\bolds{\beta}}(\hat
{\bolds{\theta}})$
is beyond the scope of this paper.

Let $\bolds\mu_0=(\mu_{0}(s_1), \ldots, \mu_{0}(s_n))^{\prime}$ and $\bolds\epsilon=(\epsilon(s_1), \ldots, \epsilon(s_n))^{\prime}$. By
$\mathbf{M}(\bolds\theta)'\bolds\Sigma^{-1}(\bolds\theta)\mathbf
{M}(\bolds\theta)
= \bolds\Sigma^{-1}(\bolds\theta)\mathbf{M}(\bolds\theta)$, (\ref
{setupfullmodelF}) and (\ref{loglikefuntrue}), we have
%
\begin{eqnarray}
-2\ell(\bolds\theta) &=& -2 \ell_{0}(\bolds\theta)+\bolds
\mu_0'\bolds\Sigma^{-1}(\bolds\theta) \bigl(
\mathbf{I}-\mathbf{M}(\bolds\theta)\bigr)\bolds\mu_0\nonumber
\\
\label{loglikedecompose}
&&{}+2\bolds\mu_0'\bolds\Sigma^{-1}(\bolds
\theta) \bigl(\mathbf{I}-\mathbf{M}(\bolds\theta)\bigr) (\bolds\eta
+\bolds
\epsilon)
\\
\nonumber
&&{}-(\bolds\eta+\bolds\epsilon)'\bolds\Sigma^{-1}(\bolds
\theta)\mathbf{M}(\bolds\theta) (\bolds\eta+\bolds\epsilon),
\end{eqnarray}
where with
$h(\bolds\theta)=(\bolds\eta+\bolds\epsilon)'\bolds\Sigma^{-1}(\bolds
\theta)(\bolds\eta+\bolds\epsilon)
- \operatorname{tr}(\bolds\Sigma^{-1}(\bolds\theta)\bolds\Sigma(\bolds\theta_0))$,
%
\begin{eqnarray}\label{logpdfe+e}
&& \ell_{0}(\bolds\theta) \equiv-\tfrac{1}{2}\bigl\{n\log(2\pi)+
\log\det\bigl(\bolds\Sigma(\bolds\theta)\bigr) +\operatorname{tr}\bigl(\bolds
\Sigma^{-1}(\bolds\theta)\bolds\Sigma(\bolds\theta_0)
\bigr)+h(\bolds\theta)\bigr\},
\end{eqnarray}
is the log-density function for $\bolds\eta+\bolds\epsilon$.
As will be seen in Section~\ref{sec2},
the contribution of the time trend
to $-2\ell(\bolds\theta)$ is mainly made by
%
\begin{eqnarray}\label{bound}
&& \bolds\mu'_0\bolds\Sigma^{-1}(\bolds\theta)
\bigl(\mathbf{I}-\mathbf{M}(\bolds\theta)\bigr)\bolds\mu_0- (\bolds
\eta+\bolds\epsilon)^{\prime}\bolds\Sigma^{-1}(\bolds\theta)
\mathbf{M}(\bolds\theta) (\bolds\eta+\bolds\epsilon).
\end{eqnarray}
The first term above, vanishing when (\ref{correct}) holds true, is due to
model misspecification, and the second term, having an order of
magnitude $O_{p}(p_{n})$ uniformly over $\Theta$ (see Lemma~\ref{lemmaprojectionsup}), is related to model complexity.
We therefore introduce
%
\begin{eqnarray}\label{eqrisk}
&& R(\Theta) =\max\Bigl\{\sup_{\bolds\theta\in\Theta}\bolds\mu'_0
\bolds\Sigma^{-1}(\bolds\theta) \bigl(\mathbf{I}-\mathbf{M}(\bolds\theta)
\bigr)\bolds\mu_0, p_{n}\Bigr\}, 
\end{eqnarray}
as a uniform bound for (\ref{bound}) over $\Theta$.
Let $(\hat{\theta}_1,\hat{\theta}_2,\hat{\theta}_3)'=\hat{\bolds\theta}$.
The growing rates of $D$
needed for $\hat{\theta}_i, i=1,2, 3$, to achieve consistency are given
in the next theorem
in terms of the order of magnitude of $R(\Theta)$.
It provides a preliminary answer to the question of whether
the covariance structures of $\bolds\eta$ and $\bolds\epsilon$
can be learnt from data under possible model misspecification.



\begin{thm}\label{thmconsistency}
Suppose
%
\begin{eqnarray}
\label{assumptionofmisselected}
&& R(\Theta)=O_{p}\bigl(n^\xi\bigr),
\end{eqnarray}
for some $\xi\in[0,1)$. Then, for $\delta\in[0,1)$,
%
\begin{eqnarray}\label{lemmatheta1consistency}
\hat\theta_1 & = & \theta_{0,1} + o_p(1)\qquad
\mbox{if }0\leq\xi<1,
\\
\label{lemmatheta2consistency}
\hat\theta_2 &= & \theta_{0,2} + o_p(1)\qquad
\mbox{if }0\leq\xi<(1+\delta)/2,
\\
\label{lemmatheta3consistency}
\hat\theta_3 &= & \theta_{0,3} + o_p(1)\qquad
\mbox{if }0\leq\xi<\delta.
\end{eqnarray}
\end{thm}

Theorem~\ref{thmconsistency} shows that as long as (\ref
{assumptionofmisselected}) holds true,
$\hat{\theta}_1$ is a consistent estimator of $\theta_{0,1}$,
regardless of the value of $\delta$.
In contrast, in order for $\hat{\theta}_{2}$
and $\hat{\theta}_{3}$ to achieve consistency,
one would require $0\leq\xi<(1+\delta)/2$
and $0\leq\xi<\delta$, respectively.
In fact, these two constraints cannot be weakened because
we provide counterexamples in Section~\ref{sec3} illustrating
that $\hat{\theta}_{3}$ is no longer consistent when
$\xi=\delta$, and both
$\hat{\theta}_{2}$ and $\hat{\theta}_{3}$
fail to achieve consistency if
$\xi=(1+\delta)/2$.
It is worth mentioning that
$\ell(\bolds\theta)$ is highly convoluted
due to the involvement of regression terms, making it difficult to establish
consistency of $\hat{\bolds{\theta}}$. Our strategy is to decompose the
nonstochastic part of
$-2\ell(\bolds\theta)$ into several layers
whose first three leading orders are
$n_{1}\equiv n$, $n_{2} \equiv n^{(1+\delta)/2}$ and $n_{3} \equiv
n^{\delta}$, respectively,
and express the remainder stochastic part
as the sum of
$h(\bolds\theta)$
and two other terms that can be uniformly expressed as
$O_{p}(R(\Theta))$ and $o_{p}( n^{\delta})$;
see (\ref{lemmaconvergea.s.0}). One distinctive
characteristic of these nonstochastic layers is that the coefficient
associated with the $i$th ($1\leq i\leq3$) leading layer only
depends on $\theta_1, \dots, \theta_{i}$. When (\ref
{assumptionofmisselected}) is assumed, this hierarchical layer structure
together with some uniform bounds established for the second moments
of $h(\bolds{\theta})$ enables us to derive the consistency of
$\hat{\bolds{\theta}}$ in the order of $\hat{\theta}_1$,
$\hat{\theta}_2$ and $\hat{\theta}_3$ by focusing on one layer and
one parameter at a time.
Let $\operatorname{tr}(\mathbf{A})$ denote the trace of a
matrix $\mathbf{A}$. As shown in the proof of Theorem~\ref
{thmconsistency}, the
uniform bounds for $h(\bolds{\theta})$ are first expressed in terms the
supremums of $\operatorname{tr}\{(\partial^{m}
\bolds\Sigma^{-1}(\bolds\theta)\bolds\Sigma(\bolds\theta_{0})/\partial
\theta_{j_{1}}\cdots\theta_{j_{m}})^{2}\}, 1\leq m \leq3,
j_{1}<\cdots<j_{m} \in\{1, 2, 3\}$, or other
similar trace terms such as those given in
\eqref{pfthm01eq01}. These expressions are obtained using the idea that the
sup-norms of a sufficiently smooth function can be bounded above by
suitable integral norms, as suggested in Lai \cite{Lai}, Chan and Ing
\cite{Chan} and Chan, Huang and Ing \cite{Chan2013}. We then carefully calculate
the orders of magnitude of the aforementioned traces, yielding uniform
bounds in terms of $n, n^{(1+\delta)/2}$ or $n^{\delta}$. Note that
Dahlhaus \cite{Dahlhaus} has applied the chaining lemma (see Pollard \cite{Pollard}) to
obtain uniform probability bounds for some quadratic forms of a
discrete time long-memory process. However, since no rates have been
reported in his bounds, his approach may not be directly
applicable here.

Whereas Theorem~\ref{thmconsistency} has demonstrated the
performance of $\hat{\bolds\theta}$ from the perspective of
consistency, the questions of what are the convergence rates of and
whether there are central limit theorems (CLTs) for
$\hat{\theta}_{i}, i=1,\ldots, 3$, still remain unanswered. The next
section is devoted to these questions. In particular, it is shown in
Theorem~\ref{thmctlunderincorrectmodel} that for
$n_i\rightarrow\infty$, $1\leq i\leq3$, $\hat{\theta}_i- \theta_{0,
i}=O_{p}(\max\{n^{\xi}n^{-1}_{i}, n^{-1/2}_{i}\})$ if
$n^{\xi}=o(n_{i})$, and $n_{i}^{1/2}(\hat{\theta}_i- \theta_{0, i})$
has a limiting normal distribution if $n^{\xi}=o(n^{1/2}_{i})$.
Since the time trend is involved, our proof of Theorem~\ref{thmctlunderincorrectmodel} is somewhat nonstandard. We
first obtain the initial convergence rates of $\hat{\bolds{\theta}}$
using the standard Taylor expansion and an argument similar but
subtler than the one used in the proof of Theorem~\ref{thmconsistency}. Using these initial rates, we can improve the
convergence results through the same argument. We then repeat this
iterative procedure until the final convergence results are
established.

The rest of this article is organized as follows.
In Section~\ref{sec2}, we begin by establishing
the CLT for $\hat{\theta}_{i}, i=1,\ldots, 3$
in situations where $p_{n}$
is fixed and the regression model is correctly specified (namely, (\ref
{correct}) is true); see Theorem~\ref{thmctlundertruemodel}.
We subsequently drop these two restrictions and report in Theorem~\ref
{thmctlunderincorrectmodel} the most general convergence results of
this paper.
In Section~\ref{sec3}, we provide two counterexamples showing that the results
obtained in Theorem~\ref{thmconsistency}
are difficult to improve.
The proofs of all theorems and corollaries in the first three sections
are given in Section~\ref{sectionproofoftheoremsandcorollaries}.
The proofs of the auxiliary lemmas used in Section~\ref{sectionproofoftheoremsandcorollaries} are
provided in the supplementary material (Chang, Huang and Ing \cite{Chang3}) in
light of space constraint.
Before leaving this section, we remark that
although our results are derived under the Gaussianity of $\{\eta(t)\}$
and $\{\epsilon(t)\}$,
similar results can be obtained when either $\{\eta(t)\}$ or $\{\epsilon
(t)\}$ is not (but pretended to be) Gaussian,
provided some fourth moment information is available.
On the other hand, while
we allow the time trend to be misspecified,
we preclude a misspecified covariance model.
The interested reader is referred to
Xiu \cite{X} for some asymptotic results on the ML estimators
when the covariance model considered in Stein \cite{Stein1990}
or A\"{\i}t-Sahalia, Mykland and Zhang \cite{AMZ}
is misspecified.

\section{Central limit theorems and rates of convergence}\label{sec2}

In this section, we begin with
establishing the asymptotic normality of $\hat{\theta}_{i}, 1\leq i
\leq3$, in
situations where the regression model is correctly
specified and $p_{n}$ is fixed.

\begin{thm}
\label{thmctlundertruemodel}
Assume that (\ref{correct}) holds and $p_{n}$
is a fixed nonnegative integer. (Note that these assumptions
yield $\xi=0$ in \eqref{assumptionofmisselected}.) Then for $\delta\in[0,1)$,
%
\begin{eqnarray}\label{thmtheta1consistency}
n^{1/2}(\hat\theta_1-\theta_{0,1}) & \displaystyle\mathop{\rightarrow}^{d} & N\bigl(0, 2\theta_{0,1}^2
\bigr),
\\
\label{thmtheta2consistency}
n^{(1+\delta)/4} (\hat{\theta}_2-\theta_{0,2} )& \displaystyle\mathop{
\rightarrow}^{d} & N \bigl(0, 2^{5/2}\theta_{0,1}^{1/2}
\theta_{0,2}^{3/2} \bigr),
\end{eqnarray}
and for $\delta\in(0,1)$,
%
\begin{eqnarray}
\label{thmtheta3consistency}
&& n^{\delta/2}(\hat\theta_3-\theta_{0,3}) \displaystyle\mathop{
\rightarrow}^{d} N(0, 2\theta_{0,3}).
\end{eqnarray}
\end{thm}

One of the easiest ways to understand
Theorem~\ref{thmctlundertruemodel} is to link the result to the Fisher
information matrix.
Straightforward calculations show that
under the assumption of Theorem~\ref{thmctlundertruemodel},
the diagonal elements of the Fisher information matrix
evaluated at $\bolds{\theta}=\bolds{\theta}_{0}$
are given by
%
\begin{eqnarray}
-\mathrm{E} \biggl(\frac{\partial^2}{\partial\theta_1^2}\ell(\bolds
{\theta}_0) \biggr)
&=& \frac{1}{2}\operatorname{tr}\bigl(\bolds{\Sigma}^{-2}(\bolds{
\theta}_{0})\bigr)+O(1),\nonumber
\\
 \label{Fisherview0}
-\mathrm{E} \biggl(\frac{\partial^2}{\partial\theta_2^2}\ell(\bolds
{\theta}_0) \biggr)
&=& \frac{1}{2\theta^{2}_{0,2}}\operatorname{tr} \bigl\{\bigl(\bolds{\Sigma
}^{-1}(\bolds{
\theta}_{0})\bolds{\Sigma}_{\eta}(\bolds{\theta}_{0})
\bigr)^{2} \bigr\}+O(1),
\\
-\mathrm{E} \biggl(\frac{\partial^2}{\partial\theta_3^2}\ell(\bolds
{\theta}_0) \biggr)
&=& \frac{1}{2}\operatorname{tr} \biggl\{ \biggl(\bolds{\Sigma}^{-1}(
\bolds{\theta}_{0}) \frac{\partial\bolds{\Sigma}(\bolds{\theta
}_{0})}{\partial\theta_{3}} \biggr)^{2} \biggr
\}+O(1) \qquad \mbox{if } 0< \delta<1,\nonumber
\end{eqnarray}
where the trace terms are solely contributed
by the log-density (log-likelihood) function
for $\bolds\eta+\bolds\epsilon$ (defined in (\ref{logpdfe+e})),
and the $O(1)$ terms, which vanish
if the time trend is known to be zero,
are related to the model complexity.
Moreover, by (\ref{Sigma1Sigmainverse2}), (\ref{DiffSigma1Sigmainverse})
and (\ref{Sigma1Sigmainverse4}),
\begin{eqnarray}
\lim_{n \to\infty}\frac{1}{2n}\operatorname{tr}\bigl(\bolds{
\Sigma}^{-2}(\bolds{\theta}_{0})\bigr)&=& \frac{1}{2\theta_{0,1}^2},
\nonumber\\
\label{Fisherview}
\lim_{n \to\infty}\frac{1}{2\theta^{2}_{0,2}n^{(1+\delta)/2}}\textrm
{tr} \bigl\{\bigl(\bolds{
\Sigma}^{-1} (\bolds{\theta}_{0})\bolds{\Sigma}_{\eta}(
\bolds{\theta}_{0})\bigr)^{2} \bigr\}& =& \frac{1}{2^{5/2}\theta
_{0,1}^{1/2}\theta_{0,2}^{3/2}},
\\
\lim_{n \to\infty} \frac{1}{2n^{\delta}}\operatorname{tr} \biggl\{ \biggl(
\bolds{\Sigma}^{-1}(\bolds{\theta}_{0}) \frac{\partial\bolds{\Sigma
}(\bolds{\theta}_{0})}{\partial\theta_{3}}
\biggr)^{2} \biggr\} &=& \frac{1}{2\theta_{0,3}}\qquad \mbox{if } 0<\delta<1.
\nonumber
\end{eqnarray}
It is interesting pointing out that the
denominator on the right-hand side of the first equation of~(\ref
{Fisherview}) coincides exactly with the limiting variance in (\ref
{thmtheta1consistency}). This is reminiscent
of a conventional asymptotic theory for the ML estimate which says that
the limiting
variance of the ML estimate is the reciprocal of the
corresponding Fisher information number.
On the other hand,
while the reciprocals of the right-hand sides of the second and third
identities of (\ref{Fisherview}) are the same as the limiting variances
in (\ref{thmtheta2consistency}) and (\ref{thmtheta3consistency}),
the divergence rates of the corresponding trace terms
$n^{(1+\delta)/2}$ and $n^{\delta}$
are much slower than $n$.
In fact,
they are
equal to the divergence rates of the second and third leading
layers of the nonstochastic part of $-2 \ell(\bolds\theta)$; see (\ref
{lemmaconvergea.s.0}).
These findings reveal that the amounts of information related to
$\theta_{0,i}$'s have different orders of magnitude, thereby leading to
different normalizing constants in the CLTs for
$\hat{\theta}_{i}$'s.


The next theorem improves Theorem~\ref{thmctlundertruemodel} by
deriving rates of convergence of $\hat{\theta}_{i},
1\leq i \leq3$, without requiring $\xi=0$ in \eqref{assumptionofmisselected}.
It further shows that CLTs for $\hat{\theta}_{i},
1\leq i \leq3$, are still possible if the model
misspecification/complexity associated with the time trend has an
order of magnitude smaller than $n^{1/2}, n^{(1+\delta)/4}$ and
$n^{\delta/2}$, respectively.

\begin{thm}
\label{thmctlunderincorrectmodel}
Suppose that (\ref{assumptionofmisselected})
is true.
Then for $\delta\in[0,1)$,
\begin{eqnarray*}
\hat\theta_1- \theta_{0,1} &=&\cases{
O_p\bigl(n^{-1/2}\bigr);& \quad$\mbox{if }\xi<1/2$,
\vspace*{3pt}\cr
O_p\bigl(n^{-(1-\xi)}\bigr); &\quad$\mbox{if }1/2\leq\xi<1$,}
%
\\
\hat\theta_2 - \theta_{0,2} &=& \cases{
O_p\bigl(n^{-(1+\delta)/4}\bigr);& \quad$\mbox{if }\xi<(1+\delta)/4$,
\vspace*{3pt}\cr
O_p\bigl(n^{-\{(1+\delta)/2-\xi\}}\bigr);&\quad$\mbox{if }(1+\delta)/4\leq\xi
<(1+\delta)/2$,
}
%
\end{eqnarray*}
and for $\delta\in(0,1)$,
\begin{eqnarray*}
\hat\theta_3-\theta_{0,3} &=& \cases{
O_p\bigl(n^{-\delta/2}\bigr);&\quad$\mbox{if }\xi<\delta/2$,
\vspace*{3pt}\cr
O_p\bigl(n^{-(\delta-\xi)}\bigr);&\quad$\mbox{if }\delta/2\leq\xi<
\delta$.
}
%
\end{eqnarray*}
In addition, for $\delta\in[0,1)$,
\begin{eqnarray*}
n^{1/2} (\hat\theta_1-\theta_{0,1} ) & \displaystyle\mathop{
\rightarrow}^{d} & N\bigl(0,2\theta_{0,1}^2\bigr);
\qquad\mbox{if }\xi<1/2,
\\
n^{(1+\delta)/4} (\hat\theta_2-\theta_{0,2} )& \displaystyle\mathop{
\rightarrow}^{d} & N\bigl(0,2^{5/2}\theta_{0,1}^{1/2}
\theta_{0,2}^{3/2}\bigr);\qquad \mbox{if }\xi< (1+\delta)/4,
\end{eqnarray*}
and for $\delta\in(0,1)$,
\begin{eqnarray*}
&& n^{\delta/2} (\hat\theta_3-\theta_{0,3} ) \displaystyle\mathop{
\rightarrow}^{d}  N(0,2\theta_{0,3});\qquad \mbox{if }\xi<\delta/2.
\end{eqnarray*}
\end{thm}

Recall that $n_{1}=n$, $n_{2}=n^{(1+\delta)/2}$ and
$n_{3}=n^{\delta}$. It is shown in (\ref{Fisherview0}) and
(\ref{Fisherview}) that, ignoring the constant, the amount of
information regarding $\theta_{0,i}$ contained in
$\bolds\eta+\bolds\epsilon$ is $n_{i}$, $1\leq i \leq3$. On the other
hand, as will become clear later, $n^{\xi}$ can be used to measure
the amount of information contaminated by model
misspecification/complexity (again ignoring the constant).
Therefore, the first part of Theorem~\ref{thmctlunderincorrectmodel}
delivers nothing more than the simple idea that
%
\begin{eqnarray}
&&\!\!\!\mbox{Rate of convergence of }\hat{\theta}_{i}\nonumber
\\
\label{mainidea}
&&\!\!\!\quad=\max\biggl\{ \frac{\mathrm{Amount} \ \mathrm{of} \  \mathrm{information}
\ \mathrm{contaminated} \ \mathrm{by} \ \mathrm{model} \  \mathrm{misspecification}\mbox{/}\mathrm{complexity}}{
\mathrm{Amount} \  \mathrm{of} \  \mathrm{information} \  \mathrm{regarding} \  \theta_{0,i} \
\mathrm{contained} \ \mathrm{in}  \ \bolds\eta+ \bolds\epsilon},\hspace*{6pt}\quad\quad
\\
\nonumber
&&\!\!\!\qquad{}\frac{1}{(\mathrm{Amount}  \ \mathrm{of}  \ \mathrm{information} \  \mathrm{regarding} \ \theta
_{0,i}\   \mathrm{contained} \ \mathrm{in} \ \bolds\eta+\bolds\epsilon
)^{1/2}} \biggr\},
\end{eqnarray}
provided that
%
\begin{eqnarray}
&& \mbox{Amount of information contaminated by model
misspecification/complexity}
\nonumber
\\[-8pt]
\label{maincondition}
\\[-8pt]
\nonumber
&&\quad<\mbox{Amount of information regarding } \theta_{0,i}
\mbox{ contained in } \bolds\eta+ \bolds\epsilon.
\end{eqnarray}
Note that the second term on the right-hand side of (\ref{mainidea})
is the best rate one can expect
when the time trend is known to be zero.
The second part of
Theorem~\ref{thmctlunderincorrectmodel}
further indicates that
the CLTs for $\hat{\theta}_{i}$'s in
Theorem~\ref{thmctlundertruemodel}
carry over to situations
where (\ref{maincondition}) holds with the right-hand side
replaced by its square root.
To the best of our knowledge, this is one of the most general CLTs
established for $\hat{\theta}_{i}$'s.
In the following,
we present two specific examples illustrating
how the asymptotic behavior of $\hat{\theta}_{i}$'s is affected
by the interaction between $\xi$ and $\delta$.
In the first example, the model misspecification
yields $R(\Theta)=O(n^{\delta})$, and hence
$\xi=\delta$.
According to Theorem~\ref{thmctlunderincorrectmodel},
the CLTs for $\hat{\theta}_1$ and $\hat{\theta}_2$ hold for a certain
range of $\delta$.
\begin{coro}
\label{coropolymle0}
Consider the intercept-only model of (\ref{setupfullmodel}) with
$p_n=0$. Suppose that
$\mu_0(s)=\beta_{0,0}+\beta_{0,1}n^{-\delta}s$,
where $\beta_{0,0}$ and $\beta_{0,1}$ are nonzero constants.
Then for $\delta\in[0, 1)$,
%
\begin{eqnarray}
\label{cor1xi=theta}
R(\Theta) &=& O\bigl(n^{\delta}\bigr),
\\
n^{1/2}(\hat\theta_1-\theta_{0,1})& \displaystyle\mathop{
\rightarrow}^{d} & N\bigl(0, 2\theta_{0,1}^2
\bigr);\qquad \delta\in[0,1/2),
\nonumber
\\[-8pt]
\label{cor1theta1consistency}
\\[-8pt]
\nonumber
n^{1-\delta}(\hat\theta_1- \theta_{0,1}) &= &
O_p(1); \qquad\delta\in[1/2,1),
\\
n^{(1+\delta)/4} (\hat{\theta}_2-\theta_{0,2} ) & \displaystyle\mathop{
\rightarrow}^{d} & N \bigl(0, 2^{5/2}\theta_{0,1}^{1/2}
\theta_{0,2}^{3/2} \bigr); \qquad \delta\in[0,1/3),
\nonumber
\\[-8pt]
\label{cor1theta2consistency}
\\[-8pt]
\nonumber
n^{(1-\delta)/2} (\hat\theta_2 -\theta_{0,2} ) &=&
O_p(1); \qquad\delta\in[1/3,1).
\end{eqnarray}
\end{coro}

We remark that the scaling factor
$n^{-\delta}$ is introduced for the linear term, $x_1(s)=n^{-\delta}s$,
so that
$ \frac{1}{n^\delta}\int_0^{n^\delta}(x_1(s)-\bar{x}_1)^2
\,ds$ does not\vspace*{1pt} depend on $n$, where
$\bar{x}_1= \frac{1}{n^\delta}\int_0^{n^\delta}x_1(s)\,ds$.
The model misspecification in the next example results in $R(\Theta
)=O_{p}(n^{(1+\delta)/2})$,
yielding $\xi=(1+\delta)/2$.
Therefore, $\hat{\theta}_1$ is guaranteed to be consistent in view of
Theorem~\ref{thmctlunderincorrectmodel}.

\begin{coro}\label{coroexp0}
Consider the same setup as in Corollary~\ref{coropolymle0} except that
$\mu_0(s) = \beta_{0,0}+\beta_{0,1}x(s)$, where $x(\cdot)$ is generated
from a zero-mean Gaussian spatial
process with covariance function
\begin{eqnarray*}
&& \operatorname{cov}\bigl(x(s),x\bigl(s'\bigr)\bigr) =
\frac{\theta_{1,2}}{\theta_{1,3}}\exp\bigl(-\theta_{1,3}\bigl|s-s'\bigr| \bigr);
\qquad s,s'\in\bigl[0,n^\delta\bigr], 
\end{eqnarray*}
for some constants $\theta_{1,2},\theta_{1,3}>0$. Then
for $\delta\in[0, 1)$,
%
\begin{eqnarray}\label{cor1xi=1+theta2}
R(\Theta)& =& O_{p}\bigl(n^{(1+\delta)/2}\bigr),
\\
\label{cor3theta1}
\hat\theta_1 &= & \theta_{0,1}+O_p
\bigl(n^{-(1-\delta)/2}\bigr).
\end{eqnarray}
\end{coro}

It is worth noting that $\hat{\theta}_{3}$ is inconsistent under the
setup of Corollary~\ref{coropolymle0}. Moreover, both
$\hat{\theta}_{2}$ and $\hat{\theta}_{3}$ are inconsistent under the
setup of Corollary~\ref{coroexp0}. These inconsistency results will
be reported in detail in the next section. Before closing this
section we remark that our theoretical results on
$\hat{\bolds{\theta}}$ can be used to make statistical inference about
the regression function. For example, when (\ref{correct}) holds and
$p_{n} \geq1$ is a fixed integer, the convergence rate of
$\hat{\bolds{\theta}}$ obtained in Theorem~\ref{thmctlundertruemodel} plays an indispensable
role in analyzing the convergence rate of the ML estimator,
$\hat{\bolds{\beta}}(\hat{\bolds{\theta}})$, of $\bolds{\beta}$.
Recently, by
making use of Theorems \ref{thmctlundertruemodel} and \ref{thmctlunderincorrectmodel},
Chang, Huang and Ing \cite{Chang2}
established the first model selection consistency result under the
mixed domain asymptotic framework. Moreover, some technical results
established in the proofs of Theorems \ref{thmctlundertruemodel} and \ref{thmctlunderincorrectmodel} have been used by
Chang, Huang and Ing \cite{Chang2} to develop a model selection consistency
result under a misspecified covariance model.

\section{Counterexamples}\label{sec3}
Using the examples constructed in Corollaries \ref{coropolymle0} and
\ref{coroexp0},
we show in this section that
the constraints $\xi<\delta$
and $\xi<(1+\delta)/2$
imposed in Theorem~\ref{thmconsistency}
for the consistency of $\hat{\theta}_{3}$
and $\hat{\theta}_{2}$, respectively,
cannot be relaxed.

\begin{coro}\label{coropolymle}
Under the setup of Corollary~\ref{coropolymle0},
%
\begin{eqnarray}\label{cor1theta3inconsistent}
\hat\theta_3 &=& \frac{12\theta_{0,2}}{12\theta_{0,2}+\beta_{0,1}^2\theta
_{0,3}}\theta_{0,3}
+o_p(1); \qquad \delta\in(0,1).
\end{eqnarray}
\end{coro}

\begin{coro}\label{coroexp}
Under the setup of Corollary~\ref{coroexp0},
%
\begin{eqnarray}\label{cor3theta2inconsistent}
\hat\theta_2 & =& \theta_{0,2} + \theta_{1,2}
\beta_{0,1}^2 + o_p(1);\qquad \delta\in[0,1),
\\
\label{cor3theta3inconsistent}
\hat\theta_3 &=& \frac{\theta_{0,2}+\beta_{0,1}^2\theta_{1,2}}{\beta
_{0,1}^2\theta_{1,2}\theta_{1,3}^{-1}
+\theta_{0,3}\theta_{0,3}^{-1}}+ o_p(1);\qquad \delta
\in(0,1).
\end{eqnarray}
\end{coro}

\begin{figure}
\begin{tabular}{@{}cc@{}}

\includegraphics{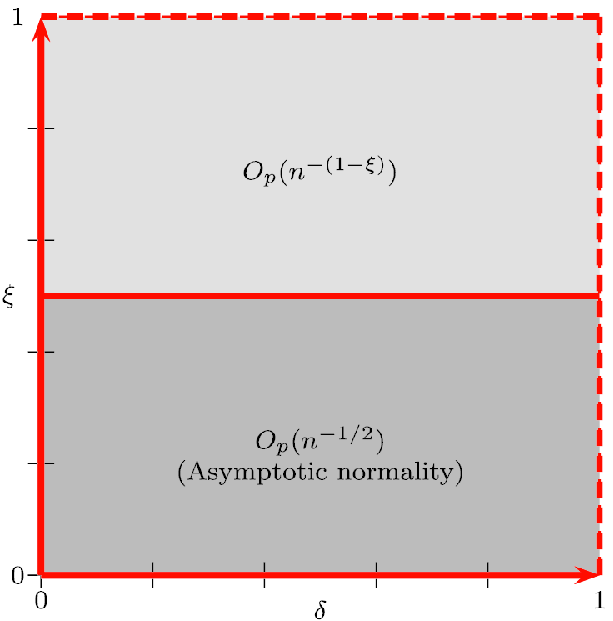}
 & \includegraphics{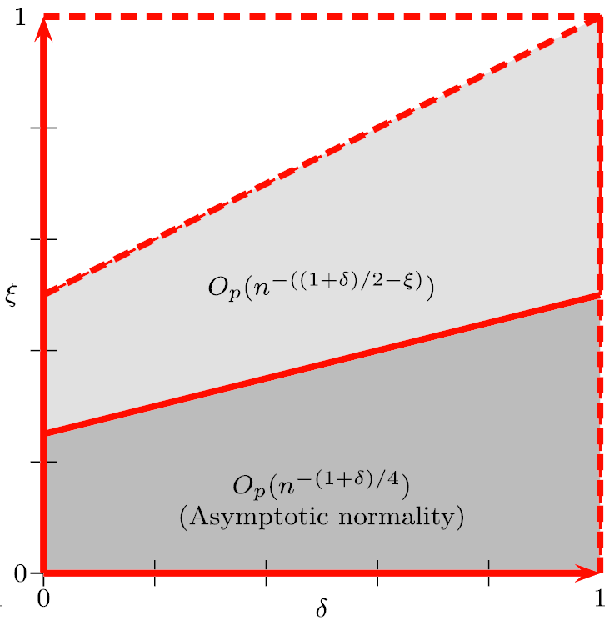}\\
\footnotesize{(a) Convergence rates of $\hat{\theta}_{1}-\theta
_{0,1}$} &
\footnotesize{(b) Convergence rates of $\hat{\theta}_{2}-\theta
_{0,2}$}\\[6pt]
\multicolumn{2}{@{}c@{}}{
\includegraphics{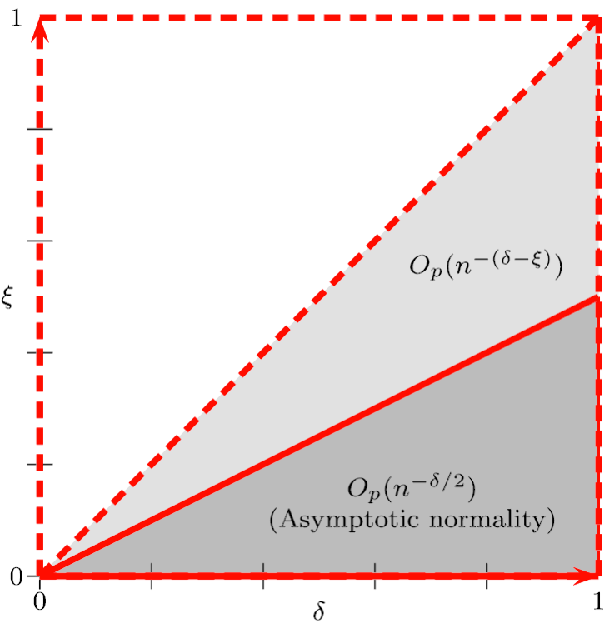}
}\\
\multicolumn{2}{@{}c@{}}{\footnotesize{(c) Convergence rates of $\hat{\theta}_{3}-\theta_{0,3}$}}
\end{tabular}
\caption{Convergence rates of $\hat{\theta}_{i}$ to $\theta_{0, i}$
with respect to $(\delta,\xi)$,
where $i=1,\ldots, 3$, $\delta$ is the growing rate of the domain and
$\xi$ satisfies $R(\Theta)=O_p(n^\xi)$.
Note that
$\hat{\theta}_{i}$ also possesses asymptotic normality
when $(\delta,\xi)$ falls in the dark gray regions,
but may fail to achieve consistency
when $(\delta,\xi)$ falls in
the white regions or
on the dash lines.
In addition, the points on the lines
between the light and dark gray area are referred to as the change
points merely in the modes of convergence but not in the convergence
rate scenario.}\vspace*{30pt}
\label{orderarea}
\end{figure}

All the above results can be illustrated by Figure~\ref{orderarea},
in which some change point behavior of $\hat{\theta}_{i}$'s (in terms
of modes of convergence)
is exhibited when $(\delta, \xi)$ runs through the region $[0,1)\times[0,1)$.

%

\section{Proofs of the theorems and corollaries}
\label{sectionproofoftheoremsandcorollaries}
In this section, we first prove the consistency of $\hat{\bolds\theta}$
in Section~\ref{sec41}.
The proofs of CLTs for $\hat{\bolds\theta}$
with and without the restrictions of
correct specification and fixed dimension
on the time trend model are given in Sections~\ref{sec42} and \ref{sec43}, respectively.
The proofs of Corollaries \ref{coropolymle0} and \ref{coropolymle}
and those of Corollaries \ref{coroexp0} and
\ref{coroexp}
are provided in Sections~\ref{sec44} and \ref{sec45}, respectively.

\subsection{Proof of Theorem \protect\ref{thmconsistency}}\label{sec41}

To prove Theorem~\ref{thmconsistency}, we need a series of
auxiliary lemmas, Lemmas \ref{lemmaardecomposition}--\ref
{lemmauniform}. Lemma~\ref{lemmaardecomposition} gives a modified
Cholesky decomposition for
$\bolds{\Sigma}^{-1}(\bolds\theta)$, which can be used to prove Lemma~\ref{propositionmaxandminboundsofeigenvalues}, asserting that
the eigenvalues of
$\bolds{\Sigma}^{-1}(\bolds\theta)\bolds{\Sigma}(\bolds\theta_{0})$ are
uniformly bounded above and below. Lemmas \ref{lem02Ttraces} and
\ref{lem03dcompseGetaG} provide the orders of magnitude of the
Cholesky factors of $\bolds{\Sigma}^{-1}(\bolds\theta)$ and the products
of $\bolds{\Sigma}_{\eta}(\bolds\theta)$ and these factors. Based on
Lemmas \ref{propositionmaxandminboundsofeigenvalues}--\ref
{lem03dcompseGetaG}, Lemma~\ref{lemmaexpequations} establishes
asymptotic expressions for the key components
of the nonstochastic part of $-2\ell(\bolds\theta)$, and Lemma~\ref{lemmaexpequations1} provides the orders of magnitude of
$\bolds{\Sigma}^{-1}(\bolds\theta)\,\partial
\bolds{\Sigma}(\bolds\theta)/\partial\theta_{i}; i=1,2,3$.
Lemmas \ref{propositionmaxandminboundsofeigenvalues}
and \ref{lemmaexpequations1} can be used in conjunction with
Lemma~\ref{lemmauniform}, which provides uniform bounds for
quadratic forms in i.i.d. random variables, to analyze the
asymptotic behavior of $h(\bolds\theta)$; see (\ref{lemmah0result}).
Lemmas \ref{lemmaprojectionsup} and \ref{lemmamisselectedmodel} explore
the effects of the time trend model on $-2
\ell(\bolds\theta)$.

%
\begin{lem}
\label{lemmaardecomposition}
Let $\bolds\Sigma(\bolds\theta)$ be given by (\ref{Sigma}) with
$\theta_1\geq0$, $\theta_2>0$ and $\theta_3>0$. Then
%
\begin{eqnarray}\label{Sigmadecomp1}
&& \bolds\Sigma^{-1}(\bolds\theta) = \mathbf{G}_n(\bolds
\theta)'\mathbf{T}_n^{-1}(\bolds\theta)
\mathbf{G}_n(\bolds\theta),
\end{eqnarray}
where
\begin{eqnarray*}
\mathbf{G}_n(\bolds\theta) & \equiv& \pmatrix{
1
& 0 & 0 & \cdots& 0
\cr
-\rho_n & 1 & 0 & \ddots& \vdots
\cr
0 & -\rho_n & 1 & \ddots& 0
\cr
\vdots& \ddots& \ddots& \ddots& 0
\cr
0 & \cdots& 0 & -\rho_n & 1
}_{n\times n},
\\
\mathbf{T}_n(\bolds\theta) &=& \mathbf{D}_n(
\bolds\theta)+\theta_1\mathbf{G}_n(\bolds\theta)
\mathbf{G}_n(\bolds\theta)', 
\end{eqnarray*}
$\rho_n=\exp(-\theta_3 n^{-(1-\delta)})$, and
\begin{eqnarray*}
&& \mathbf{D}_n(\bolds\theta)\equiv\frac{\theta_2}{\theta_3}\pmatrix{
1 & 0 & \cdots& 0
\cr
0 & 1-\rho_n^2 & \ddots& \vdots
\cr
\vdots& \ddots& \ddots& 0
\cr
0 & \cdots& 0 & 1-\rho_n^2
}_{n\times n}. 
\end{eqnarray*}
\end{lem}

\begin{lem}
\label{propositionmaxandminboundsofeigenvalues}
Let $\lambda_{\max}(\mathbf{A})$
and $\lambda_{\min}(\mathbf{A})$
denote the maximum and minimum eigenvalues
of the matrix $\mathbf{A}$. For $\bolds\Sigma(\bolds\theta)$ given by
(\ref{Sigma}),
suppose that $\Theta\subset(0,\infty)^3$ is compact.
Then,
%
\begin{eqnarray}
0 & < & \liminf_{n\rightarrow\infty} \inf_{\bolds\theta\in\Theta}
\lambda_{\min} \bigl( \bolds\Sigma^{-1/2}(\bolds\theta)\bolds
\Sigma(\bolds\theta_0) \bolds\Sigma^{-1/2}(\bolds\theta)
\bigr)
\nonumber
\\[-8pt]
\label{propositionboundeigen}
\\[-8pt]
\nonumber
& \leq& \limsup_{n\rightarrow\infty}\sup_{\bolds\theta\in\Theta}
\lambda_{\max} \bigl(\bolds\Sigma^{-1/2}(\bolds\theta) \bolds
\Sigma(\bolds\theta_0)\bolds\Sigma^{-1/2}(\bolds\theta)
\bigr)< \infty.
\end{eqnarray}
\end{lem}

\begin{lem}
\label{lem02Ttraces}
Under the setup of Lemma~\ref{lemmaardecomposition},
for any $\bolds\theta\in\Theta\in(0,\infty)^3$, where $\Theta$ is
compact, and $\delta\in[0,1)$,
the following equation holds uniformly over $\Theta$:
%
\begin{eqnarray}\label{traceT2}
&& \operatorname{tr}\bigl(\mathbf{T}_n^{-2}(\bolds\theta)
\bigr) = \frac{n^{(5-3\delta)/2}}{2^{7/2}\theta_1^{1/2}\theta_2^{3/2}}
+ o\bigl(n^{(5-3\delta)/2}\bigr).
\end{eqnarray}
\end{lem}

\begin{lem}
\label{lem03dcompseGetaG}
Under the setup of Lemma~\ref{lem02Ttraces},
for any $\bolds\theta\in\Theta$,\vspace*{-3pt}
%
\begin{eqnarray}
&& \mathbf{G}_n(\bolds\theta)\bolds\Sigma_\eta(\bolds
\theta_0)\mathbf{G}_n(\bolds\theta)'\nonumber
\\
\label{lem03eq01Geta0G}
&&\quad= \frac{\theta_{0,2}\rho_n}{\theta_{0,3}\rho_{0,n}}\bigl(1-\rho
_{0,n}^2\bigr)\mathbf{I}
+ \biggl(1-\frac{\rho_n}{\rho_{0,n}} \biggr) (1-\rho_n\rho_{0,n})
\bolds\Sigma_\eta(\bolds\theta_0)
\\
&&\qquad{}+\frac{\theta_{0,2}}{\theta_{0,3}} \biggl\{ \biggl(1-\frac{\rho
_n}{\rho_{0,n}} \biggr) \bigl(
\mathbf{v}_0\mathbf{e}_1'+
\mathbf{e}_1\mathbf{v}_0'\bigr)+
\rho_n^2\mathbf{e}_1\mathbf{e}_1'
\biggr\},\nonumber
\end{eqnarray}
where $\mathbf{e}_1=(1,0,\ldots,0)'$,
$\mathbf{v}_0=(1,\rho_{0,n},\ldots,\rho_{0,n}^{n-1})$ and $\rho
_{0,n}=\exp(-\theta_{0,3}n^{-(1-\delta)})$.
In addition, for any $\delta\in[0,1)$,\vspace*{-6pt}
%
\begin{eqnarray}\label{lem03eq02vTv}
\sup_{\bolds\theta\in\Theta}\mathbf{v}_0'
\mathbf{T}_n^{-1}(\bolds\theta)\mathbf{v}_0 &=&
O\bigl(n^{2(1-\delta)}\bigr),
\\
\label{lem03eq03vTe}
\sup_{\bolds\theta\in\Theta}\mathbf{v}_0'
\mathbf{T}_n^{-1}(\bolds\theta)\mathbf{e}_1 &=&
O\bigl(n^{1-\delta}\bigr),
\\
\label{lem03eq04eTe}
\sup_{\bolds\theta\in\Theta}\mathbf{e}_1'
\mathbf{T}_n^{-1}(\bolds\theta)\mathbf{e}_1 &=&
O(1).  
\end{eqnarray}
Furthermore,\vspace*{-6pt} for any $\delta\in(0,1)$,
%
\begin{eqnarray}\label{lem03eq06traceTetaTeta}
\sup_{\bolds\theta\in\Theta}\operatorname{tr} \bigl( \bigl(\mathbf{T}_n^{-1}(
\bolds\theta)\bolds\Sigma_\eta(\bolds\theta) \bigr)^2
\bigr) &=& \frac{1}{4\theta_3^3}n^{4-3\delta} + o\bigl(n^{4-3\delta}\bigr).
\end{eqnarray}
\end{lem}

\begin{lem}
\label{lemmaexpequations}
Under the setup of Lemma~\ref{lem02Ttraces}, the following
equations hold uniformly over\vspace*{-3pt} $\Theta$:
%
\begin{eqnarray}
\log\bigl(\det\bigl(\bolds\Sigma(\bolds\theta)\bigr)\bigr) &= & n\log
\theta_1 + \biggl(\frac{2\theta_2}{\theta_1} \biggr)^{1/2}n^{(1+\delta)/2}
- \biggl(\frac{\theta_2}{\theta_1}+\theta_3 \biggr)n^{\delta}
\nonumber
\\[-9pt]
\label{logdetSigma}
\\[-9pt]
\nonumber
&&{}-\frac{1-\delta}{2}\log n + o\bigl(n^{\delta}\bigr)+ O(1),
\\
\operatorname{tr} \bigl(\bolds\Sigma(\bolds\theta_0)\bolds
\Sigma^{-1}(\bolds\theta) \bigr)& =& \frac{\theta_{0,1}}{\theta_1}n -
\frac{\theta_{0,1}}{2\theta_1} \biggl(\frac{2\theta_2}{\theta_1} \biggr
)^{1/2}n^{(1+\delta)/2}\nonumber
\\
\label{Sigma1Sigmainverse3}
&&{}+\frac{\theta_{0,2}}{(2\theta_1\theta_2)^{1/2}}n^{(1+\delta)/2} +\frac
{\theta_{0,2}(\theta_3^2-\theta_{0,3}^2)}{2\theta_2\theta_{0,3}}n^\delta
\\
&&{}+ o\bigl(n^\delta\bigr)+O(1).\nonumber
\end{eqnarray}
\end{lem}

\begin{lem}
\label{lemmaexpequations1}
Under the setup of Lemma~\ref{lem02Ttraces}, the following
equations hold uniformly over\vspace*{-3pt} $\Theta$:
%
\begin{eqnarray}\label{Sigma1Sigmainverse2}
\operatorname{tr} \bigl( \bigl(\bolds\Sigma_\eta(\bolds\theta) \bolds
\Sigma^{-1}(\bolds\theta) \bigr)^2 \bigr)& =& \biggl(
\frac{\theta_2}{8\theta_1} \biggr)^{1/2}n^{(1+\delta)/2} +o\bigl
(n^{(1+\delta)/2}
\bigr),
\\[-2pt]
\operatorname{tr} \bigl(\bolds\Sigma_\eta(\bolds\theta_0)
\bolds\Sigma^{-1}(\bolds\theta) \bigr)& =& \frac{\theta_{0,2}}{(2\theta
_1\theta_2)^{1/2}}n^{(1+\delta)/2}
+\frac{\theta_{0,2}(\theta_3^2-\theta_{0,3}^2)}{2\theta_2\theta
_{0,3}}n^\delta
\nonumber
\\[-9pt]
\label{Sigma1Sigmainverse}
\\[-9pt]
\nonumber
&&{}+ o\bigl(n^\delta\bigr)+ O(1),
\\
\label{DiffSigma1Sigmainverse}
\operatorname{tr} \biggl( \biggl(\bolds\Sigma^{-1}(\bolds\theta)
\frac{\partial}{\partial\theta_3} \bolds\Sigma(\bolds\theta) \biggr)^2
\biggr) &= &
\frac{1}{\theta_3}n^\delta+ o\bigl(n^\delta\bigr).
\end{eqnarray}
\end{lem}

\begin{rem}
As will be shown later, (\ref{propositionboundeigen}), (\ref
{Sigma1Sigmainverse2})
and (\ref{Sigma1Sigmainverse})
can be used to derive bounds for $\operatorname{tr}((\bolds\Sigma^{-1}(\bolds
\theta)\,
\partial\bolds\Sigma(\bolds\theta)/\partial\theta_2
)^2)$ and $\operatorname{tr}((\bolds\Sigma^{-1}(\bolds\theta)
\,\partial\bolds\Sigma(\bolds\theta)/\partial\theta_1
)^2)$. These bounds, together with (\ref{DiffSigma1Sigmainverse}),
play important roles in establishing the consistency of $\hat{\theta}_{1}$.
\end{rem}

\begin{lem}\label{lemmaprojectionsup}
Let $\mathbf{X}$ be full rank a.s. Then under the setup of
Lemma~\ref{lem02Ttraces},
%
\begin{eqnarray} \label{finitevariable}
&& \sup_{\bolds\theta\in\Theta} \bigl\{ (\bolds\eta+\bolds\epsilon)'
\bolds\Sigma^{-1}(\bolds\theta)\mathbf{M}(\bolds\theta) (\bolds\eta+
\bolds\epsilon) \bigr\}= O_p(p_n),
\end{eqnarray}
where $\mathbf{M}(\bolds\theta)$ is defined in (\ref
{projectionmatrix}).
\end{lem}

\begin{lem}
\label{lemmamisselectedmodel}
Under the setup up of Lemma~\ref{lem02Ttraces}, let
$\mathbf{X}$ be full rank a.s. Suppose that for some $\xi\geq0$,
\begin{eqnarray*}
&& \sup_{\bolds\theta\in\Theta} \bigl\{\bolds\mu'_0\bolds
\Sigma^{-1}(\bolds\theta) \bigl(\mathbf{I}-\mathbf{M}(\bolds\theta)\bigr)
\bolds\mu_0 \bigr\}= O_{p}\bigl(n^\xi\bigr).
\end{eqnarray*}
Then
%
\begin{eqnarray}\label{lemmamis2}
&& \sup_{\bolds\theta\in\Theta} \bigl\{\bolds\mu'_0\bolds
\Sigma^{-1}(\bolds\theta) \bigl(\mathbf{I}-\mathbf{M}(\bolds\theta)\bigr)
(\bolds\eta+\bolds\epsilon) \bigr\} = o_p\bigl(n^\xi
\bigr).
\end{eqnarray}
\end{lem}

Before introducing Lemma~\ref{lemmauniform}, we need some
notation. For $1\leq m\leq r<\infty$, define $\mathbf{J}(m,
r)=\{(j_{1}, \ldots, j_{m}): j_{1}<\cdots<j_{m}, j_{i}\in
\{1, \ldots, r\}, 1\leq i \leq m \}$. Let $g(\bolds\xi)$ be a function
of $\bolds\xi=(\xi_1,\ldots,\xi_r)'\in\mathbb{R}^{r}$.
For $\mathbf{j}=(j_{1}, \ldots, j_{m}) \in\mathbf{J}(m, r)$, define
$\mathbf{D}_{\mathbf{j}} g(\bolds\xi)=
\partial^{m} g(\bolds\xi)/\partial\xi_{j_{1}}, \ldots, \partial\xi_{j_{m}}$.
Denote by $B_{\tau}(\bolds\lambda)$ the $r$-dimensional closed
ball centered at $\bolds\lambda=(\lambda_{1}, \ldots,
\lambda_{r})^{\prime}$ with radius $0<\tau<\infty$. For $ {\mathbf
j}\in\mathbf{J}(m, r)$, define the $m$-dimensional sphere:
\begin{eqnarray*}
B_{\tau}(\bolds\lambda, {\mathbf j}) &=& \bigl\{(\xi_{j_1},
\ldots,\xi_{j_m}): (\lambda_{1},\ldots,\lambda_{j_1-1},
\xi_{j_1},\lambda_{j_1+1},\ldots,\lambda_{j_2-1},
\\
&&{}\xi_{j_2},\lambda_{j_2+1},\ldots,
\lambda_{j_m-1}, \xi_{j_m},\lambda_{j_m+1},\ldots,
\lambda_{r})\in B_{\tau}(\bolds\lambda)\bigr\}.
\end{eqnarray*}

\begin{lem}
\label{lemmauniform}
Assume that $w_{1}, \ldots, w_{n}$ are i.i.d. random variables with
$\mathrm{E}(w_{1})=0, \mathrm{E}(w^{2}_{1})=1$ and
$\mathrm{E}(w^{4}_{1})<\infty$. Let
$\mathbf{A}(\bolds\xi)=[a_{i,j}(\bolds\xi)]_{1\leq i, j\leq n}$ be an
$n \times n$ matrix whose $(i,j)$th component is $a_{i,j}(\bolds\xi)$,
a function of $\bolds\xi$ with a continuous partial derivative
$\mathbf{D}_{\mathbf{j}} a_{i,j}(\bolds\xi)$ on
$B_{\tau}(\bolds\lambda)$, for $\mathbf{j} \in
\mathbf{J}(m,r)$. Define
$q_1(\bolds\xi)=\mathbf{w}^{\prime}\mathbf{A}(\bolds\xi)\mathbf
{w}-\operatorname{tr}(\mathbf{A}(\bolds\xi))$,
where $\mathbf{w}=(w_{1}, \ldots, w_{n})'$. Then for
$\bolds\xi\in B_{\tau}(\bolds\lambda)$, there exists a constant $C>0$
such that
%
\begin{eqnarray}
&& \mathrm{E} \Bigl(\sup_{\bolds\xi\in B_{\tau}(\bolds\lambda)} \bigl
(q_1(\bolds
\xi)-q_1(\bolds\lambda)\bigr)^{2} \Bigr)
\nonumber
\\[-8pt]
\label{lemmauniformbound}
\\[-8pt]
\nonumber
&&\quad\leq  C \sum_{m=1}^{r}\sum
_{\mathbf{j} \in
\mathbf{J}(m,r)}\operatorname{vol}^{2}\bigl(B_{\tau}(
\bolds\lambda, {\mathbf j})\bigr) \sup_{\bolds\xi\in B_{\tau}(\bolds
\lambda)} \operatorname{var}\bigl(
\mathbf{D}_{\mathbf{j}}q_1(\bolds\xi)\bigr),
\end{eqnarray}
where $\operatorname{vol}(\Theta)$ denotes the volume of $\Theta$.
\end{lem}

First, we prove (\ref{lemmatheta1consistency}).
By (\ref{assumptionofmisselected}), (\ref{logdetSigma}), (\ref
{Sigma1Sigmainverse3}), (\ref{finitevariable}) and (\ref{lemmamis2}),
it follows that
%
\begin{eqnarray}
-2\ell(\bolds\theta) &=& n\log(2\pi) - \frac{1-\delta}{2}\log{n} + \biggl
(\log
\theta_1 + \frac{\theta_{0,1}}{\theta_1} \biggr)n
\nonumber\\
&&{} + \biggl(\frac{2\theta_2}{\theta_1} \biggr)^{1/2} \biggl(1-\frac
{\theta_{0,1}}{2\theta_1}
+\frac{\theta_{0,2}}{2\theta_2} \biggr)n^{(1+\delta)/2}
\nonumber
\\[-8pt]
\label{lemmaconvergea.s.0}
\\[-8pt]
\nonumber
&&{}- \biggl\{\frac{\theta_2}{\theta_1} +\theta_3-\frac{\theta
_{0,2}(\theta_3^2 - \theta_{0,3}^2)}{2\theta_2\theta_{0,3}} \biggr\}
n^\delta
\\
&&{}+h(\bolds\theta) + O_p\bigl(n^\xi\bigr) +
o_p\bigl(n^\delta\bigr),\nonumber
\end{eqnarray}
uniformly in $\Theta$, where
$h(\bolds\theta)=(\bolds\eta+\bolds\epsilon)'\bolds\Sigma^{-1}(\bolds
\theta)(\bolds\eta+\bolds\epsilon)
- \operatorname{tr}(\bolds\Sigma^{-1}(\bolds\theta)\bolds\Sigma(\bolds\theta_0))$.
Hence, (\ref{lemmatheta1consistency}) is
ensured by for any $\varepsilon>0$,
%
\begin{eqnarray}\label{lemmatheta1consistency0}
&& P \Bigl(\inf_{\bolds\theta\in\Theta_1(\epsilon)}\bigl\{-2\ell(\bolds
\theta)+2\ell(\bolds
\theta_0)\bigr\}>0 \Bigr)\rightarrow1,
\end{eqnarray}
as $n\rightarrow\infty$, where $\Theta_1(\epsilon)
= \{\bolds\theta\in\Theta:|\theta_1-\theta_{0,1}|>\varepsilon\}$.
Since by (\ref{lemmaconvergea.s.0}),
\begin{eqnarray*}
\inf_{\bolds\theta\in\Theta_1(\varepsilon)}\bigl\{-2\ell(\bolds\theta)
+ 2\ell(\bolds
\theta_0)\bigr\} & \geq& \inf_{\bolds\theta\in\Theta_1(\varepsilon)}
\biggl\{\log
\theta_1 +\frac{\theta_{0,1}}{\theta_1} -\log(\theta_{0,1}) - 1 \biggr
\}n
\\
&&{}-\sup_{\bolds\theta\in\Theta_1(\varepsilon)}\bigl|h(\bolds\theta) -
h(\bolds\theta_0)\bigr|
+o_p(n),
\end{eqnarray*}
and since
$ \inf_{\bolds\theta\in\Theta_1(\varepsilon)} \{\log\theta_1
+\frac{\theta_{0,1}}{\theta_1} -\log(\theta_{0,1}) - 1 \}>0$,
(\ref{lemmatheta1consistency0}) follows immediately from
%
\begin{equation}\label{lemmah0result}
\mbox{E} \Bigl(\sup_{\bolds\theta\in\Theta} \bigl|h(\bolds\theta)-h(\bolds
\theta_0) \bigr|^2 \Bigr) = O(n).
\end{equation}

Since
$h(\bolds\theta)$
is continuous on
$\Theta$ and $\Theta$ is compact, in the rest of the proof, we assume
without loss of
generality that $\Theta=B_{\tau}(\bolds\theta_{0})$, a closed ball centered
at $\bolds\theta_{0}$ with radius $\tau$ for some $0<\tau<\infty$.
By (\ref{lemmauniformbound}) with
$\mathbf{w}=\bolds\Sigma^{-1/2}(\bolds\theta_0)(\bolds\eta+\bolds
\epsilon)$ and
$\mathbf{A}(\bolds\theta)=\bolds\Sigma^{1/2}(\bolds\theta_0)\bolds\Sigma
^{-1}(\bolds\theta)\bolds\Sigma^{1/2}(\bolds\theta_0)$,
we obtain
$h(\bolds\theta)=\mathbf{w}^{\prime}\mathbf{A}(\bolds\theta)\mathbf
{w}-\operatorname{tr}(\mathbf{A}(\bolds\theta))$
and
\begin{eqnarray}
&& \mbox{E} \Bigl(  \sup_{\bolds\theta\in\Theta} \bigl|h(\bolds\theta)-h(\bolds
\theta_0) \bigr|^2 \Bigr)\nonumber
\\
&&\quad\leq  C\sup_{\bolds\theta\in\Theta} \biggl\{ \operatorname{var} \biggl(
\frac{\partial}{\partial\theta_1} h(\bolds\theta) \biggr) + \operatorname{var} \biggl(
\frac{\partial}{\partial\theta_2} h(\bolds\theta) \biggr) + \operatorname{var} \biggl(
\frac{\partial}{\partial\theta_3} h(\bolds\theta) \biggr)\nonumber
\nonumber
\\[-8pt]
\label{Ing0504a}
\\[-8pt]
\nonumber
&&\qquad{}+ \operatorname{var} \biggl(\frac{\partial^2}{\partial\theta_1\,\partial\theta
_2} h(\bolds\theta) \biggr) +\operatorname{var}
\biggl(\frac{\partial^2}{\partial\theta_1\,\partial\theta_3} h(\bolds
\theta) \biggr) +\operatorname{var} \biggl(
\frac{\partial^2}{\partial\theta_2\,\partial\theta_3} h(\bolds\theta)
\biggr)
\\
&&\qquad{}+\operatorname{var} \biggl(\frac{\partial^3}{\partial\theta_1\,\partial\theta
_2\,\partial\theta_3} h(\bolds\theta) \biggr) \biggr\},\nonumber
\end{eqnarray}
for some constant $C>0$.
By (\ref{propositionboundeigen}), (\ref{Sigma1Sigmainverse}),
%
\begin{equation}\label{Ing0504b}
\operatorname{tr}(A)\lambda_{\min}(B)\leq\operatorname{tr}(AB) \leq\operatorname{tr}(A)
\lambda_{\max}(B),
\end{equation}
for the nonnegative definite matrices $A$ and $B$,
and using $\mathbf{I}-\bolds\Sigma^{-1}(\bolds\theta_0)\Sigma_\eta
(\bolds\theta_0)=\theta_0\bolds\Sigma^{-1}(\bolds\theta_0)$
twice, we obtain
%
\begin{eqnarray}
\operatorname{tr}\bigl(\bolds\Sigma^{-2}(\bolds\theta_0)\bigr)
&= & \frac{1}{\theta_{0,1}} \bigl\{ \operatorname{tr}\bigl(\bolds\Sigma
^{-1}(\bolds
\theta_0)\bigr)-\operatorname{tr}\bigl(\bolds\Sigma^{-2} (\bolds
\theta_0)\bolds\Sigma_\eta(\bolds\theta_0)
\bigr) \bigr\}
\nonumber\\
&= & \frac{1}{\theta_{0,1}} \biggl\{ \frac{1}{\theta_{0,1}}\bigl(n -
\operatorname{tr}\bigl(
\bolds\Sigma^{-1}(\bolds\theta_0) \bolds
\Sigma_\eta(\bolds\theta_0)\bigr)\bigr)
\nonumber
\\[-8pt]
\label{Sigma1Sigmainverse4}
\\[-8pt]
\nonumber
&&{}-\operatorname{tr}\bigl(\bolds\Sigma^{-2} (\bolds\theta_0)
\bolds\Sigma_\eta(\bolds\theta_0)\bigr) \biggr\}
\\
&= & \frac{1}{\theta_{0,1}^2}n + O\bigl(n^{(1+\delta)/2}\bigr).\nonumber
\end{eqnarray}
Equations (\ref{propositionboundeigen}), (\ref{Ing0504b}) and (\ref
{Sigma1Sigmainverse4}) lead to
%
\begin{eqnarray}
\sup_{\bolds\theta\in\Theta} \operatorname{var} \biggl(\frac{\partial
}{\partial\theta_1} h(\bolds
\theta) \biggr) &=& \sup_{\bolds\theta\in\Theta} 2\operatorname{tr} \biggl(
\biggl(
\frac{\partial}{\partial\theta_1} \bolds\Sigma^{-1}(\bolds\theta)\bolds
\Sigma(\bolds
\theta_0) \biggr)^2 \biggr)
\nonumber
\\[-8pt]
\label{Ing0505a}
\\[-8pt]
\nonumber
&=& \sup_{\bolds\theta\in\Theta} 2\operatorname{tr} \bigl( \bigl(\bolds
\Sigma^{-2}(\bolds\theta) \bolds\Sigma(\bolds\theta_0)
\bigr)^2 \bigr)= O(n).
\end{eqnarray}
Similarly, (\ref{propositionboundeigen}), (\ref{Sigma1Sigmainverse2})
and (\ref{Ing0504b}) imply
%
\begin{eqnarray}
\sup_{\bolds\theta\in\Theta} \operatorname{var} \biggl(\frac{\partial
}{\partial\theta_2} h(\bolds
\theta) \biggr) &=& \sup_{\bolds\theta\in\Theta} 2\operatorname{tr} \biggl(
\biggl(
\frac{\partial}{\partial\theta_2} \bolds\Sigma^{-1}(\bolds\theta)\bolds
\Sigma(\bolds
\theta_0) \biggr)^2 \biggr)\nonumber
\\
\label{Ing0505b}
&=& \sup_{\bolds\theta\in\Theta} \frac{2}{\theta_2^2}\operatorname{tr}
\bigl(
\bigl(\bolds\Sigma^{-1}(\bolds\theta) \bolds\Sigma_\eta(\bolds
\theta)\bolds\Sigma^{-1}(\bolds\theta)\bolds\Sigma(\bolds
\theta_0) \bigr)^2 \bigr)
\\
&=& O\bigl(n^{(1+\delta)/2}\bigr).\nonumber 
\end{eqnarray}
Moreover,
by (\ref{propositionboundeigen}), (\ref{DiffSigma1Sigmainverse})
and (\ref{Ing0504b}), one gets
%
\begin{eqnarray}
\sup_{\bolds\theta\in\Theta} \operatorname{var} \biggl(\frac{\partial
}{\partial\theta_3} h(\bolds
\theta) \biggr)& =& \sup_{\bolds\theta\in\Theta} 2\operatorname{tr} \biggl(
\biggl(
\frac{\partial}{\partial\theta_3}\bolds\Sigma^{-1}(\bolds\theta) \bolds
\Sigma(\bolds
\theta_0) \biggr)^2 \biggr)\nonumber
\\
\label{Ing0505c}
&=& \sup_{\bolds\theta\in\Theta} 2\operatorname{tr} \biggl( \biggl(\bolds
\Sigma^{-1}(\bolds\theta) \biggl( \frac{\partial}{\partial\theta
_3}\bolds\Sigma(\bolds
\theta) \biggr) \bolds\Sigma^{-1}(\bolds\theta)\bolds\Sigma(\bolds
\theta_0) \biggr)^2 \biggr)
\\
&=& O\bigl(n^{\delta}\bigr). 
\nonumber
\end{eqnarray}
In a similar way, it can be shown that
%
\begin{eqnarray}\label{Ing0505d}
&& \sup_{\bolds\theta\in\Theta} \operatorname{var} \biggl(\frac{\partial
^2}{\partial\theta_1\,\partial\theta_2} h(\bolds
\theta) \biggr)=O\bigl(n^{(1+\delta)/2}\bigr),
\end{eqnarray}
and
%
\begin{eqnarray}
&& \sup_{\bolds\theta\in\Theta} \operatorname{var} \biggl(\frac{\partial
^2}{\partial\theta_1\,\partial\theta_3} h(\bolds
\theta) \biggr)+ \sup_{\bolds\theta\in\Theta} \operatorname{var} \biggl(
\frac{\partial^2}{\partial\theta_2\,\partial\theta_3} h(\bolds\theta)
\biggr)
\nonumber
\\[-8pt]
\label{Ing0505e}
\\[-8pt]
\nonumber
&&\quad{}+ \sup_{\bolds\theta\in\Theta} \operatorname{var} \biggl(\frac{\partial
^3}{\partial\theta_1\,\partial\theta_2\,\partial\theta_3} h(\bolds
\theta) \biggr)=O\bigl(n^{\delta}\bigr).
\end{eqnarray}
Consequently, (\ref{lemmah0result})
follows from
(\ref{Ing0504a})--(\ref{Ing0505e}),
and hence
(\ref{lemmatheta1consistency}) holds true.

Next, we prove (\ref{lemmatheta2consistency}), which in turn is implied by
the property that for any $\varepsilon_2>0$, there exists an
$\varepsilon_1>0$ such that
%
\begin{eqnarray}\label{theta2consistent}
&& \mathrm{P} \Bigl(\inf_{\bolds\theta\in\Theta_2(\bolds{\varepsilon
})}\bigl\{ -2\ell(\bolds\theta)+2\ell
\bigl((\theta_1,\theta_{0,2},\theta_{0,3})'
\bigr)\bigr\}>0 \Bigr)\rightarrow1,
\end{eqnarray}
as $n\rightarrow\infty$, where
$\Theta_2(\bolds{\varepsilon})=
\{\bolds\theta\in\Theta:|\theta_1-\theta_{0,1}|\leq\varepsilon_1,
|\theta_2-\theta_{0,2}|>\varepsilon_2\}$ and $\bolds\varepsilon
=(\varepsilon_1,\varepsilon_2)'$. Let
$\bolds\theta_b=(\theta_1,\theta_{0,2},\theta_{0,3})'$. Since
$\xi<(1+\delta)/2$, by (\ref{lemmaconvergea.s.0}), we have
\begin{eqnarray*}
&& \inf_{\bolds\theta\in\Theta_2(\bolds{\varepsilon})} \bigl\{-2\ell
(\bolds\theta) + 2\ell(\bolds
\theta_b)\bigr\}
\\
&&\quad\geq  \inf_{\bolds\theta\in\Theta_2(\bolds{\varepsilon})} \frac
{1}{(2\theta_1\theta_2)^{1/2}} \biggl\{ \bigl(
\theta_2^{1/2}-\theta_{0,2}^{1/2}
\bigr)^2 +\theta_2^{1/2}\bigl(
\theta_2^{1/2}-\theta_{0,2}^{1/2}\bigr)
\biggl(1-\frac{\theta_{0,1}}{\theta_1} \biggr) \biggr\}n^{(1+\delta)/2}
\\
&&\qquad {}-\sup_{\bolds\theta\in\Theta_2(\bolds{\varepsilon})}\bigl|h(\bolds\theta)
- h(\bolds\theta_b)\bigr|
+o_p\bigl(n^{(1+\delta)/2}\bigr).
\end{eqnarray*}
Therefore (\ref{theta2consistent}) is given by
%
\begin{eqnarray}\label{lemmaconvergea.s.5}
&& \mathrm{E} \Bigl(\sup_{\bolds\theta\in\Theta}\bigl|h(\bolds\theta) - h(\bolds
\theta_b)\bigr|^2 \Bigr) = O_p
\bigl(n^{(1+\delta)/2}\bigr).
\end{eqnarray}
By
(\ref{lemmauniformbound}) with
$\mathbf{w}=\bolds\Sigma^{-1/2}(\bolds\theta_0)(\bolds\eta+\bolds
\epsilon)$ and
$\mathbf{A}(\bolds\theta)=\bolds\Sigma^{1/2}(\bolds\theta_0)
\{\bolds\Sigma^{-1}(\bolds\theta)-\bolds\Sigma^{-1}(\bolds\theta_b) \}
\bolds\Sigma^{1/2}(\bolds\theta_0)$,
we obtain
$h(\bolds\theta)-h(\bolds\theta_b)=\mathbf{w}^{\prime}\mathbf{A}(\bolds
\theta)
\mathbf{w}-\operatorname{tr}(\mathbf{A}(\bolds\theta))$ and
%
\begin{eqnarray}
&& \mbox{E} \Bigl(  \sup_{\bolds\theta\in\Theta} \bigl|h(\bolds\theta)-h(\bolds
\theta_b) \bigr|^2 \Bigr)\nonumber
\\
&&\quad\leq C\sup_{\bolds\theta\in\Theta} \biggl\{ \operatorname{var} \biggl(
\frac{\partial}{\partial\theta_1} \bigl(h(\bolds\theta)-h(\bolds\theta
_b) \bigr)
\biggr) + \operatorname{var} \biggl(\frac{\partial}{\partial\theta_2} h(\bolds
\theta) \biggr) +
\operatorname{var} \biggl(\frac{\partial}{\partial\theta_3} h(\bolds\theta)
\biggr)
\nonumber
\\[-8pt]
\label{Ing0505f}
\\[-8pt]
\nonumber
&&\qquad{}+ \operatorname{var} \biggl(\frac{\partial^2}{\partial\theta_1\,\partial\theta
_2} h(\bolds\theta) \biggr) +\operatorname{var}
\biggl(\frac{\partial^2}{\partial\theta_1\,\partial\theta_3} h(\bolds
\theta) \biggr) +\operatorname{var} \biggl(
\frac{\partial^2}{\partial\theta_2\,\partial\theta_3} h(\bolds\theta)
\biggr)
\\
&&\qquad{}+\operatorname{var} \biggl(\frac{\partial^3}{\partial\theta_1\,\partial\theta
_2\,\partial\theta_3} h(\bolds\theta) \biggr) \biggr\},\nonumber
\end{eqnarray}
for some constant $C>0$.
In addition,
it follows from
(\ref{propositionboundeigen}), (\ref{Sigma1Sigmainverse2}) and (\ref
{Ing0504b}) that
%
\begin{eqnarray}
&& \sup_{\bolds\theta\in\Theta}\operatorname{var} \biggl(
\frac{\partial}{\partial\theta_1} \bigl(h(\bolds\theta)-h(\bolds\theta
_b) \bigr)
\biggr)\nonumber
\\
&&\quad= 2 \sup_{\bolds\theta\in\Theta}\operatorname{tr} \biggl( \biggl(\frac
{\partial}{\partial\theta_1}
\bigl( \bolds\Sigma^{-1}(\bolds\theta)-\bolds\Sigma^{-1}(\bolds
\theta_b) \bigr) \bolds\Sigma(\bolds\theta_0)
\biggr)^2 \biggr)\nonumber
\\
&&\quad= 2 \sup_{\bolds\theta\in\Theta} \operatorname{tr} \bigl( \bigl(\bigl(\bolds
\Sigma^{-2}(\bolds\theta)-\bolds\Sigma^{-2}(\bolds
\theta_b)\bigr) \bolds\Sigma(\bolds\theta_0)
\bigr)^2 \bigr)\nonumber
\\
&&\quad= 2 \sup_{\bolds\theta\in\Theta}\operatorname{tr} \bigl( \bolds\Sigma^{1/2}(
\bolds\theta_0)\bolds\Sigma^{-2}(\bolds\theta) \bigl(\bolds
\Sigma^2(\bolds\theta_b)-\bolds\Sigma^2(
\bolds\theta)\bigr)\bolds\Sigma^{-2}(\bolds\theta_b)\bolds
\Sigma(\bolds\theta_0)\nonumber
\\
&&\qquad {}\times\bolds\Sigma^{-2}(\bolds\theta_b) \bigl(\bolds
\Sigma^2(\bolds\theta_b)-\bolds\Sigma^2(
\bolds\theta)\bigr) \bolds\Sigma^{-2}(\bolds\theta)\bolds
\Sigma^{1/2}(\bolds\theta_0) \bigr)\nonumber
\\
&&\quad=  O \Bigl(\sup_{\bolds\theta\in\Theta}\operatorname{tr} \bigl( \bolds
\Sigma^{1/2}(\bolds\theta_0)\bolds\Sigma^{-2}(
\bolds\theta) \bigl(\bolds\Sigma^2(\bolds\theta_b)-\bolds
\Sigma^2(\bolds\theta)\bigr)\nonumber\\
&&\qquad\quad{}\times \bolds\Sigma^{-2}(\bolds
\theta_b) \bigl(\bolds\Sigma^2(\bolds
\theta_b)-\bolds\Sigma^2(\bolds\theta)\bigr) \bolds
\Sigma^{-2}(\bolds\theta) \bolds\Sigma^{1/2}(\bolds
\theta_0) \bigr) \Bigr)\nonumber
\\
&&\quad= O \Bigl(\sup_{\bolds\theta\in\Theta}\operatorname{tr} \bigl( \bolds
\Sigma^{-2}(\bolds\theta)\bolds\Sigma(\bolds\theta_0)\bolds
\Sigma^{-2}(\bolds\theta) \bigl(\bolds\Sigma^2(\bolds
\theta_b)-\bolds\Sigma^2(\bolds\theta)\bigr) \bolds
\Sigma^{-2}(\bolds\theta_b) \bigl(\bolds
\Sigma^2(\bolds\theta_b)-\bolds\Sigma^2(
\bolds\theta)\bigr) \bigr) \Bigr)\nonumber
\\
&&\quad=  O \Bigl(\sup_{\bolds\theta\in\Theta}\operatorname{tr} \bigl( \bolds
\Sigma^{-2}(\bolds\theta) \bigl(\bolds\Sigma^2(\bolds
\theta_b)-\bolds\Sigma^2(\bolds\theta)\bigr) \bolds
\Sigma^{-2}(\bolds\theta_b) \bigl(\bolds
\Sigma^2(\bolds\theta_b)-\bolds\Sigma^2(
\bolds\theta)\bigr) \bigr) \Bigr)\nonumber
\\
&&\quad= O \Bigl(\sup_{\bolds\theta\in\Theta}\operatorname{tr} \bigl( \bigl(\bolds
\Sigma^{-2}(\bolds\theta)-\bolds\Sigma^{-2}(\bolds
\theta_b)\bigr) \bigl(\bolds\Sigma^2(\bolds
\theta_b)-\bolds\Sigma^2(\bolds\theta)\bigr) \bigr)
\Bigr)\nonumber
\\
&&\quad= O \Bigl(\sup_{\bolds\theta\in\Theta}\operatorname{tr} \bigl( \bigl(\bolds
\Sigma^{-1}(\bolds\theta) \bigl(\bolds\Sigma^{-1}(\bolds\theta)-
\bolds\Sigma^{-1}(\bolds\theta_b)\bigr) +\bigl(\bolds
\Sigma^{-1}(\bolds\theta)-\bolds\Sigma^{-1}(\bolds
\theta_b)\bigr)\bolds\Sigma^{-1}(\bolds\theta_b)
\bigr)
\nonumber
\\[-8pt]
\label{pfthm01eq01}
\\[-8pt]
\nonumber
&&\qquad {}\times\bigl(\bolds\Sigma(\bolds\theta_b) \bigl(\bolds\Sigma(
\bolds\theta_b)-\bolds\Sigma(\bolds\theta)\bigr) +\bigl(\bolds\Sigma(
\bolds\theta_b)-\bolds\Sigma(\bolds\theta)\bigr)\bolds\Sigma(\bolds
\theta) \bigr) \bigr) \Bigr)
\\
&&\quad= O \Bigl(\sup_{\bolds\theta\in\Theta}\operatorname{tr} \bigl( \bolds
\Sigma^{-1}(\bolds\theta) \bigl(\bolds\Sigma^{-1}(\bolds\theta)-
\bolds\Sigma^{-1}(\bolds\theta_b)\bigr)\bolds\Sigma(\bolds
\theta_b) \bigl(\bolds\Sigma(\bolds\theta_b)-\bolds
\Sigma(\bolds\theta)\bigr) \bigr) \Bigr)\nonumber
\\
&&\qquad {}+O \Bigl(\sup_{\bolds\theta\in\Theta}\operatorname{tr} \bigl( \bolds
\Sigma^{-1}(\bolds\theta) \bigl(\bolds\Sigma^{-1}(\bolds\theta)-
\bolds\Sigma^{-1}(\bolds\theta_b)\bigr) \bigl(\bolds\Sigma(
\bolds\theta_b)-\bolds\Sigma(\bolds\theta)\bigr)\bolds\Sigma(\bolds
\theta) \bigr) \Bigr)\nonumber
\\
&&\qquad {}+O \Bigl(\sup_{\bolds\theta\in\Theta}\operatorname{tr} \bigl( \bigl(\bolds
\Sigma^{-1}(\bolds\theta)-\bolds\Sigma^{-1}(\bolds
\theta_b)\bigr)\bolds\Sigma^{-1}(\bolds\theta_b)
\bolds\Sigma(\bolds\theta_b) \bigl(\bolds\Sigma(\bolds
\theta_b)-\bolds\Sigma(\bolds\theta)\bigr) \bigr) \Bigr)\nonumber
\\
&&\qquad{}+O \Bigl(\sup_{\bolds\theta\in\Theta}\operatorname{tr} \bigl( \bigl(\bolds
\Sigma^{-1}(\bolds\theta)-\bolds\Sigma^{-1}(\bolds
\theta_b)\bigr)\bolds\Sigma^{-1}(\bolds\theta_b)
\bigl(\bolds\Sigma(\bolds\theta_b)-\bolds\Sigma(\bolds\theta)\bigr)
\bolds\Sigma(\bolds\theta) \bigr) \Bigr)\nonumber
\\
&&\quad= O \Bigl(\sup_{\bolds\theta\in\Theta}\operatorname{tr} \bigl( \bolds
\Sigma^{-2}(\bolds\theta) \bigl(\bolds\Sigma(\bolds\theta_b)-
\bolds\Sigma(\bolds\theta)\bigr)^2 \bigr) \Bigr)\nonumber
\\
&&\qquad {}+O \Bigl(\sup_{\bolds\theta\in\Theta}\operatorname{tr} \bigl( \bolds
\Sigma^{-1}(\bolds\theta) \bigl(\bolds\Sigma(\bolds\theta_b)-
\bolds\Sigma(\bolds\theta)\bigr) \bolds\Sigma^{-1}(\bolds
\theta_b) \bigl(\bolds\Sigma(\bolds\theta_b)-\bolds
\Sigma(\bolds\theta)\bigr) \bigr) \Bigr)\nonumber
\\
&&\qquad{}+O \Bigl(\sup_{\bolds\theta\in\Theta}\operatorname{tr} \bigl( \bolds
\Sigma^{-1}(\bolds\theta) \bigl(\bolds\Sigma(\bolds\theta_b)-
\bolds\Sigma(\bolds\theta)\bigr) \bolds\Sigma^{-1}(\bolds
\theta_b) \bigl(\bolds\Sigma(\bolds\theta_b)-\bolds
\Sigma(\bolds\theta)\bigr) \bigr) \Bigr)\nonumber
\\
&&\qquad {}+O \Bigl(\sup_{\bolds\theta\in\Theta}\operatorname{tr} \bigl( \bolds
\Sigma^{-2}(\bolds\theta_b) \bigl(\bolds\Sigma(\bolds
\theta_b)-\bolds\Sigma(\bolds\theta)\bigr)^2 \bigr)
\Bigr)\nonumber
\\
&&\quad= O \Bigl(\sup_{\bolds\theta\in\Theta}\operatorname{tr} \bigl( \bolds
\Sigma^{-2}(\bolds\theta) \bigl(\bolds\Sigma_\eta(\bolds
\theta_b)-\bolds\Sigma_\eta(\bolds\theta)
\bigr)^2 \bigr) \Bigr)\nonumber
\\
&&\quad= O\bigl(n^{(1+\delta)/2}\bigr).\nonumber
\end{eqnarray}
Combining (\ref{Ing0505f})
and (\ref{pfthm01eq01}), with
(\ref{Ing0505b})--(\ref{Ing0505e}), yields (\ref{lemmaconvergea.s.5}),
and hence
(\ref{lemmatheta2consistency}) is established.

Finally, we prove (\ref{lemmatheta3consistency}). It
suffices to show that for any $\varepsilon_3>0$, there exist
$\varepsilon_1,\varepsilon_2>0$ such that
%
\begin{eqnarray}\label{theta3consistent}
&& \mathrm{P} \Bigl(\inf_{\bolds\theta\in\Theta_3(\bolds{\varepsilon
})}\bigl\{-2\ell(\bolds\theta) +2\ell
\bigl((\theta_1,\theta_2,\theta_{0,3})'
\bigr)\bigr\}> 0 \Bigr)\rightarrow1,
\end{eqnarray}
as $n\rightarrow\infty$, where
$\Theta_3(\bolds{\varepsilon})=\{\bolds\theta\in\Theta:|\theta_1-\theta
_{0,1}|\leq\varepsilon_1,
|\theta_2-\theta_{0,2}|\leq\varepsilon_2,|\theta_{3}-\theta
_{0,3}|>\varepsilon_3\}$
and $\bolds{\varepsilon}=(\varepsilon_1,\varepsilon_2,\varepsilon_3)'$. Let
$\bolds\theta_c=(\theta_1,\theta_2,\theta_{0,3})'$. Since $\xi<\delta$,
by (\ref{lemmaconvergea.s.0}), we have
\begin{eqnarray*}
\inf_{\bolds\theta\in\Theta_3(\bolds{\varepsilon})} \bigl\{-2\ell(\bolds
\theta) +2\ell(\bolds
\theta_c)\bigr\} & \geq& \inf_{\bolds\theta\in\Theta_3(\bolds{\varepsilon
})} \biggl\{
\frac{\theta_{0,2}(\theta_3 - \theta_{0,3})^2}{2\theta_{0,3}\theta_2}
-(\theta_3-\theta_{0,3}) \biggl(1-
\frac{\theta_{0,2}}{\theta_2} \biggr) \biggr\}n^\delta
\\
&&{}-\sup_{\bolds\theta\in\Theta_3(\bolds{\varepsilon})}\bigl|h(\bolds\theta
)-h(\bolds{\theta}_c)\bigr|
+ o_p\bigl(n^\delta\bigr).
\end{eqnarray*}
Therefore, it suffices for (\ref{theta3consistent}) to
show that
%
\begin{eqnarray}\label{lemmaconsistencyhsmallo}
&& \mathrm{E} \Bigl(\sup_{\bolds\theta\in\Theta} \bigl|h(\bolds\theta)-h(\bolds
\theta_c) \bigr|^2 \Bigr) = O\bigl(n^\delta\bigr).
\end{eqnarray}
By (\ref{lemmauniformbound}) with
$\mathbf{w}=\bolds\Sigma^{-1/2}(\bolds\theta_0)(\bolds\eta+\bolds
\epsilon)$ and
$\mathbf{A}(\bolds\theta)=\bolds\Sigma^{1/2}(\bolds\theta_0)
\{\bolds\Sigma^{-1}(\bolds\theta)-\bolds\Sigma^{-1}(\bolds\theta_c) \}
\bolds\Sigma^{1/2}(\bolds\theta_0)$,
we obtain
$h(\bolds\theta)-h(\bolds\theta_c)=\mathbf{w}^{\prime}\mathbf{A}(\bolds
\theta)
\mathbf{w}-\operatorname{tr}(\mathbf{A}(\bolds\theta))$ and
%
\begin{eqnarray}
&& \mathrm{E} \Bigl(  \sup_{\bolds\theta\in\Theta} \bigl|h(\bolds\theta
)-h(\bolds
\theta_c) \bigr|^2 \Bigr)
\nonumber\\
&&\quad\leq C\sup_{\bolds\theta\in\Theta} \biggl\{ \operatorname{var} \biggl(
\frac{\partial}{\partial\theta_1} \bigl(h(\bolds\theta)-h(\bolds\theta
_c) \bigr)
\biggr) +\operatorname{var} \biggl(\frac{\partial}{\partial\theta_2} \bigl
(h(\bolds\theta)-h(\bolds
\theta_c) \bigr) \biggr) +\operatorname{var} \biggl(\frac{\partial}{\partial
\theta_3}h(
\bolds\theta) \biggr)\qquad
\nonumber
\\[-8pt]
\label{Ing0505g}
\\[-8pt]
\nonumber
&&\qquad{}+\operatorname{var} \biggl(\frac{\partial^2}{\partial\theta_1\,\partial\theta
_2} \bigl(h(\bolds\theta)-h(\bolds
\theta_c) \bigr) \biggr) +\operatorname{var} \biggl(\frac{\partial^2}{\partial
\theta_1\,\partial\theta_3}h(
\bolds\theta) \biggr) +\operatorname{var} \biggl(\frac{\partial^2}{\partial
\theta_2\,\partial\theta_3}h(\bolds\theta)
\biggr)
\\
\nonumber
&&\qquad {}+\operatorname{var} \biggl(\frac{\partial^3}{\partial\theta_1\,\partial\theta
_2\,\partial\theta_3}h(\bolds\theta) \biggr) \biggr\},
\end{eqnarray}
for some constant $C>0$.
In view of (\ref{Ing0505g}),
(\ref{Ing0505c}) and (\ref{Ing0505e}),
(\ref{lemmaconsistencyhsmallo}) is guaranteed by
%
\begin{eqnarray}\label{lemmaconsistencyhderivatetheta1}
\sup_{\bolds\theta\in\Theta} \operatorname{var} \biggl(\frac{\partial
}{\partial\theta_1} \bigl(h(
\bolds\theta)-h(\bolds\theta_c) \bigr) \biggr)& =& O
\bigl(n^\delta\bigr),
\\
\label{lemmaconsistencyhderivatetheta2}
\sup_{\bolds\theta\in\Theta} \operatorname{var} \biggl(\frac{\partial
}{\partial\theta_2} \bigl(h(
\bolds\theta)-h(\bolds\theta_c) \bigr) \biggr) &=& O
\bigl(n^\delta\bigr),
\\[-2pt]
\label{lemmaconsistencyhderivatetheta12}
\sup_{\bolds\theta\in\Theta} \operatorname{var} \biggl(\frac{\partial
^2}{\partial\theta_1\,\partial\theta_2} \bigl(h(
\bolds\theta)-h(\bolds\theta_c) \bigr) \biggr) &=& O
\bigl(n^\delta\bigr).
\end{eqnarray}
In what follows, we only focus on the proof of (\ref
{lemmaconsistencyhderivatetheta1}) since the proofs of
(\ref{lemmaconsistencyhderivatetheta2}) and
(\ref{lemmaconsistencyhderivatetheta12}) are similar.
Note first that
by an argument similar to that used to prove
\eqref{pfthm01eq01},\vspace*{-6pt} one obtains
%
\begin{eqnarray}
&&\sup_{\bolds\theta\in\Theta} \operatorname{var} \biggl(\frac{\partial
}{\partial\theta_1} \bigl(h(
\bolds\theta)-h(\bolds\theta_c) \bigr) \biggr)\nonumber
\\[-2pt]
&&\quad= 2 \sup_{\bolds\theta\in\Theta}\operatorname{tr} \biggl( \biggl(\frac
{\partial}{\partial\theta_1}
\bigl( \bolds\Sigma^{-1}(\bolds\theta)-\bolds\Sigma^{-1}(\bolds
\theta_c) \bigr) \bolds\Sigma(\bolds\theta_0)
\biggr)^2 \biggr)
\nonumber
\\[-9pt]
\label{pfthm01eq02}
\\[-9pt]
\nonumber
&&\quad= O \Bigl( \sup_{\bolds\theta\in\Theta}\operatorname{tr} \bigl( \bigl(\bolds
\Sigma^{-1}(\bolds\theta) \bigl(\bolds\Sigma(\bolds\theta_c)-
\bolds\Sigma(\bolds\theta)\bigr) \bigr)^2 \bigr) \Bigr)
\\[-2pt]
&&\quad= O \Bigl(\sup_{\bolds\theta\in\Theta}\operatorname{tr} \bigl( \bigl(
\mathbf{T}_n^{-1}(\bolds\theta) \mathbf{G}_n(
\bolds\theta) \bigl(\bolds\Sigma_{\eta}(\bolds\theta_c)-
\bolds\Sigma_{\eta}(\bolds\theta) \bigr) \mathbf{G}_n(\bolds
\theta)' \bigr)^2 \bigr) \Bigr).\nonumber
\end{eqnarray}
In addition, (\ref{lem03eq01Geta0G})
and some algebraic manipulations\vspace*{-3pt} yield
%
\begin{eqnarray}
&&\mathbf{G}_n(\bolds\theta) \bigl(\bolds\Sigma_{\eta}(
\bolds\theta_c)-\bolds\Sigma_{\eta}(\bolds\theta)\bigr)
\mathbf{G}_n(\bolds\theta)'\nonumber
\\[-2pt]
\label{pfthm01eq03}
&&\quad= \biggl(\frac{\theta_2\rho_n}{\theta_{0,3}\rho_{0,n}}\bigl(1-\rho
_{0,n}^2\bigr) -
\frac{\theta_2}{\theta_3}\bigl(1-\rho_n^2\bigr) \biggr)
\mathbf{I} + \biggl(1-\frac{\rho_n}{\rho_{0,n}} \biggr) (1-\rho_n
\rho_{0,n})\bolds\Sigma_\eta(\bolds\theta_c)
\\[-2pt]
&&\qquad {}+\frac{\theta_2}{\theta_{0,3}} \biggl(1-\frac{\rho_n}{\rho_{0,n}}
\biggr) \bigl(
\mathbf{v}_0\mathbf{e}_1'+
\mathbf{e}_1\mathbf{v}_0'\bigr) +
\theta_2 \biggl(\frac{1}{\theta_{0,3}}-\frac{1}{\theta_3} \biggr)
\rho_n^2\mathbf{e}_1\mathbf{e}_1',
\nonumber
\end{eqnarray}
where $\rho_{0,n}=\exp(-\theta_{0,3} n^{-(1-\delta)})$,\vspace*{-3pt} and
%
\begin{eqnarray}\label{pfprop01eq011minusrho}
&& 1-\rho_n^k\rho_{0,n}^\ell= (k
\theta_3+\ell\theta_{0,3})n^{-(1-\delta)} + O
\bigl(n^{-2(1-\delta)}\bigr); \qquad k,\ell\in\mathbb{Z},
\end{eqnarray}
uniformly in $\Theta$.
Moreover, by (\ref{traceT2}), (\ref{lem03eq02vTv})--(\ref
{lem03eq06traceTetaTeta})
and\vspace*{-3pt}
\begin{eqnarray*}
&& \limsup_{n\rightarrow\infty}\sup_{\bolds\theta\in\Theta}
\lambda_{\max}\bigl( \bolds\Sigma_\eta^{-1}(\bolds
\theta)\bolds\Sigma_\eta(\bolds\theta_c)\bigr)<\infty,
\end{eqnarray*}
which can be shown using an argument similar to that
used to prove (B.2) in the supplementary document (Chang, Huang
and Ing \cite{Chang3}),
we\vspace*{-3pt} have
\begin{eqnarray}
\sup_{\bolds\theta\in\Theta}n^{-4(1-\delta)} \operatorname{tr} \bigl(
\mathbf{T}_n^{-2}(\bolds\theta) \bigr) &=& O
\bigl(n^\delta\bigr),
\nonumber\\[-2pt]
\sup_{\bolds\theta\in\Theta} n^{-4(1-\delta)} \operatorname{tr} \bigl(\bigl(
\mathbf{T}_n^{-1}(\bolds\theta)\bolds\Sigma_\eta(
\bolds\theta_c)\bigr)^2 \bigr) &=& O\bigl(n^\delta
\bigr),\nonumber
\nonumber
\\[-9pt]
\label{Ing0505i}
\\[-9pt]
\nonumber
\sup_{\bolds\theta\in\Theta}n^{-2(1-\delta)} \operatorname{tr} \bigl(\bigl(
\mathbf{T}_n^{-1}(\bolds\theta) \bigl(\mathbf{v}_0
\mathbf{e}_1'+\mathbf{e}_1
\mathbf{v}_0'\bigr)\bigr)^2 \bigr) &=& O(1),
\\[-2pt]
\sup_{\bolds\theta\in\Theta} \operatorname{tr} \bigl(\bigl(\mathbf{T}_n^{-1}(
\bolds\theta)\mathbf{e}_1\mathbf{e}_1'
\bigr)^2 \bigr) &=& O(1).\nonumber
\end{eqnarray}
Combining (\ref{pfthm01eq02})--(\ref{Ing0505i}) leads to
(\ref{lemmaconsistencyhderivatetheta1})
and hence (\ref{lemmaconsistencyhsmallo}). This completes the proof of
(\ref{lemmatheta3consistency}).

\subsection{Proof of Theorem
\protect\ref{thmctlundertruemodel}}\label{sec42}
To prove Theorem~\ref{thmctlundertruemodel}, we need two
additional lemmas, Lemmas \ref{lemmataylorfirstorder}--\ref
{lemmataylorsecondorder}, which provide the
orders of magnitude of
$\partial\ell(\hat{\bolds\theta})/\partial\theta_i$ and
$\partial^2\ell(\hat{\bolds\theta})/\partial\theta_i^2$; $i=1,2,3$,
when the convergence rate of $\hat{\bolds\theta}$ is given.
On the contrary, using the orders of the
magnitude of $\partial\ell(\hat{\bolds\theta})/\partial\theta_i$ and
$\partial^2\ell(\hat{\bolds\theta})/\partial\theta_i^2$, $i=1,2,3$,
one can also derive
the convergence rate of
$\hat{\bolds\theta}$; see
(\ref{thmtheta1taylor})--(\ref{thmtheta3taylor}).
As a result, the convergence rate of
$\hat{\bolds\theta}$ can be sequentially improved via
an initial convergence rate and applying this argument repeatedly.

\begin{lem}
\label{lemmataylorfirstorder}
Under the setup of Lemma~\ref{lemmamisselectedmodel}, define for $k=1,2,3$,
\begin{eqnarray*}
&& g_k(\bolds\theta) = -\frac{\partial}{\partial\theta_k}2\ell(\bolds
\theta),
\end{eqnarray*}
where $\ell(\bolds\theta)$ is given by (\ref
{loglikefuntrue}). Let $\hat{\bolds\theta}
=(\hat\theta_1,\hat\theta_2,\hat\theta_3)'$ be an estimate of
$\bolds\theta$ with $\hat\theta_1 = \theta_{0,1}+O_p(n^{-r_1})$,
$\hat\theta_2=\theta_{0,2}+ O_p(n^{-r_2})$ and
$\hat\theta_3=\theta_{0,3}+O_p(n^{-r_3})$ for some constants
$r_1\in[0,1/2]$, $r_2\in[0,(1+\delta)/4]$ and $r_3\in[0,\delta/2]$;
$\delta\in[0,1)$. Then for any $\delta\in[0,1)$,
%
\begin{eqnarray}
g_1\bigl((\theta_{0,1},\hat\theta_2,\hat
\theta_3)'\bigr) &=& O_p\bigl(n^{1/2}
\bigr) + O_p\bigl(n^{(1+\delta)/2-r_2}\bigr)+ O_p
\bigl(n^{\delta-r_3}\bigr)
\nonumber
\\[-8pt]
\label{thmg1order}
\\[-8pt]
\nonumber
&&{}+ O_p\bigl(n^\xi\bigr) + O(1),
\\
g_2\bigl((\hat\theta_1,\theta_{0,2},\hat
\theta_3)'\bigr) &=& O_p\bigl(n^{(1+\delta)/4}
\bigr) + O_p\bigl(n^{(1+\delta)/2-r_1}\bigr)+ O_p
\bigl(n^{\delta-r_3}\bigr)
\nonumber
\\[-8pt]
\label{thmg2order}
\\[-8pt]
\nonumber
&&{}+ O_p\bigl(n^\xi\bigr) + O(1),
\end{eqnarray}
and for $\delta\in(0,1)$,
%
\begin{eqnarray}
g_3\bigl((\hat\theta_1,\hat\theta_2,
\theta_{0,3})'\bigr) &=& O_p
\bigl(n^{\delta/2}\bigr) + O_p\bigl(n^{\delta-r_1}\bigr) +
O_p\bigl(n^{\delta-r_2}\bigr)
\nonumber
\\[-8pt]
\label{thmg3order}
\\[-8pt]
\nonumber
&&{}+ O_p\bigl(n^\xi\bigr) + O(1).
\end{eqnarray}
In addition, for any $\delta\in[0,1)$, if $\xi<1/2$ and
$r_2\geq\delta/2$,
%
\begin{eqnarray}\label{thmg1asymptoticnormality}
n^{-1/2}g_1\bigl((\theta_{0,1},\hat
\theta_2,\hat\theta_3)'\bigr) & \displaystyle\mathop{
\rightarrow}^{d}& N\bigl(0,2\theta_{0,1}^{-2}\bigr);
\end{eqnarray}
if $\xi<(1+\delta)/4$, $r_1>(1+\delta)/4$ and $r_3>-(1-3\delta)/4$,
%
\begin{eqnarray}\label{thmg2asymptoticnormality}
&& n^{-(1+\delta)/4} g_2\bigl((\hat\theta_1,
\theta_{0,2},\hat\theta_3)'\bigr) \displaystyle\mathop{
\rightarrow}^{d} N\bigl(0,2^{-1/2}\theta_{0,1}^{-1/2}
\theta_{0,2}^{-3/2}\bigr).
\end{eqnarray}
Furthermore, for any $\delta\in(0,1)$, if $\xi<\delta/2$,
$r_1>\delta/2$ and $r_2>\delta/2$,
%
\begin{eqnarray}\label{thmg3asymptoticnormality}
&& n^{-\delta/2}g_3\bigl((\hat\theta_1,\hat
\theta_2,\theta_{0,3})'\bigr) \displaystyle\mathop{
\rightarrow}^{d} N\bigl(0,2\theta_{0,3}^{-1}\bigr).
\end{eqnarray}
\end{lem}

\begin{lem}
\label{lemmataylorsecondorder}
Under the setup of Lemma~\ref{lemmamisselectedmodel}, let
%
\begin{eqnarray}\label{thmgii}
&& g_{kk}(\bolds\theta) = -\frac{\partial^2}{\partial\theta_k^2}2\ell
(\bolds\theta);\qquad
k=1,2,3.
\end{eqnarray}
Let $\hat{\bolds\theta}
=(\hat\theta_1,\hat\theta_2,\hat\theta_3)'$ be an estimate of $\bolds
\theta$.
Suppose that $\hat\theta_1=\theta_{0,1}+o_p(1)$. Then for $\delta\in[0,1)$,
there exists a constant $\theta_{0,1}^*>0$ satisfying
$|\theta_{0,1}^*-\hat\theta_1|\leq|\theta_{0,1}-\hat\theta_1|$ such that
%
\begin{eqnarray}\label{thmg11order}
&& g_{11}\bigl(\bigl(\theta_{0,1}^*,\hat\theta_2,
\hat\theta_3\bigr)'\bigr) = \frac{n}{\theta_{0,1}^2}
+o_p(n).
\end{eqnarray}
In addition, suppose that $\hat\theta_1=\theta_{0,1}+o_p(1)$
and $\hat\theta_2=\theta_{0,2}+o_p(1)$,
then for $\delta\in[0,1)$, there exists a constant $\theta_{0,2}^*>0$ satisfying
$|\theta_{0,2}^*-\hat\theta_2|\leq|\theta_{0,2}-\hat\theta_2|$ such that
%
\begin{eqnarray}\label{thmg22order}
&& g_{22}\bigl(\bigl(\hat\theta_1,\theta_{0,2}^*,
\hat\theta_3\bigr)'\bigr) = \frac{n^{(1+\delta)/2}}{2^{3/2}\theta
_{0,1}^{1/2}\theta_{0,2}^{3/2}} +
O_p\bigl(n^\xi\bigr)+o_p
\bigl(n^{(1+\delta)/2}\bigr).
\end{eqnarray}
Furthermore, suppose that $\hat{\bolds\theta}=\bolds\theta
_0+o_p(1)$, then for $\delta\in(0,1)$,
there exists a constant $\theta_{0,3}^*>0$ satisfying
$|\theta_{0,3}^*-\hat\theta_3|\leq|\theta_{0,3}-\hat\theta_3|$ such that
%
\begin{eqnarray}\label{thmg33order}
g_{33}\bigl(\bigl(\hat\theta_1,\hat\theta_2,
\theta_{0,3}^*\bigr)'\bigr) &=& \frac{n^\delta}{\theta_{0,3}}
+O_p\bigl(n^\xi\bigr)+o_p
\bigl(n^\delta\bigr).
\end{eqnarray}
\end{lem}

We shall prove (\ref{thmtheta1consistency})--(\ref
{thmtheta3consistency}) by iteratively applying (\ref
{thmg1order})--(\ref{thmg33order}).
For the first iteration, we show that
%
\begin{eqnarray}\label{thmtheta1iterate1}
\hat\theta_1 - \theta_{0,1}& =& O_p
\bigl(n^{-(1-\delta)/2}\bigr)\qquad \mbox{if }\delta\in[0,1),
\\
n^{(1+\delta)/4} (\hat{\theta}_2-\theta_{0,2} ) &\displaystyle\mathop{
\rightarrow}^{d} & N \bigl(0, 2^{5/2}\theta_{0,1}^{1/2}
\theta_{0,2}^{3/2} \bigr)\qquad \mbox{if }\delta\in[0,1/3),
\nonumber
\\[-8pt]
\label{thmtheta2iterat1}
\\[-8pt]
\nonumber
\hat\theta_2 - \theta_{0,2} &=& O_p
\bigl(n^{-(1-\delta)/2}\bigr)\qquad \mbox{if }\delta\in[1/3,1),
\\
n^{\delta/2} (\hat{\theta}_3-\theta_{0,3} )& \displaystyle\mathop{
\rightarrow}^{d} & N(0,2\theta_{0,3})\qquad \mbox{if }\delta
\in(0,1/2),
\nonumber
\\[-8pt]
\label{thmtheta3iterate1}
\\[-8pt]
\nonumber
\hat{\theta}_3-\theta_{0,3} &= & O_p
\bigl(n^{-(1-\delta)/2}\bigr)\qquad \mbox{if }\delta\in[1/2,1).
\end{eqnarray}

\begin{pf*}{Proof of (\protect\ref{thmtheta1iterate1})}
Taking\vspace*{-3pt} the Taylor expansion of
$g_1(\hat{\bolds\theta})$ at $\hat{\bolds\theta}_a=(\theta_{0,1},
\hat\theta_2,\hat\theta_3)'$ yields
%
\begin{equation}\label{thmtheta1taylor}
0 = g_1(\hat{\bolds\theta}) = g_1(\hat{\bolds
\theta}_a) + g_{11}\bigl(\hat{\bolds\theta}_a^*
\bigr) (\hat{\theta}_1-\theta_{0,1}),
\end{equation}
where $\hat{\bolds\theta}_a^*=(\theta_{0,1}^*,\hat\theta_2,\hat\theta
_3)'$ satisfies
$|\theta_{0,1}^*-\hat\theta_1|\leq|\theta_{0,1}-\hat\theta_1|$.
Therefore, for (\ref{thmtheta1iterate1}) to hold,
it suffices to show that\vspace*{-3pt}
\begin{eqnarray*}
g_1(\hat{\bolds\theta}_a) &=& O_p
\bigl(n^{(1+\delta)/2}\bigr),
\\[-2pt]
g_{11}\bigl(\hat{\bolds\theta_a^*}\bigr)& =&
\frac{n}{\theta_{0,1}^2} + o_p(n),
\end{eqnarray*}
where the first equation follows from (\ref{lemmatheta2consistency}) and
(\ref{thmg1order}) with $r_2=0$, and the second
one is given by (\ref{lemmatheta1consistency}) and
(\ref{thmg11order}).
\end{pf*}

\begin{pf*}{Proof of (\protect\ref{thmtheta2iterat1})}
Let
$\hat{\bolds\theta}_b=(\hat\theta_1,\theta_{0,2},\hat\theta_3)'$.
Taking the Taylor expansion of
$g_2(\hat{\bolds\theta})$ at $\hat{\bolds\theta}_b=(\hat\theta_1,\theta
_{0,2},\hat\theta_3)'$ yields
%
\begin{equation}\label{thmtheta2taylor}
0 = g_2(\hat{\bolds\theta}) = g_2(\hat{\bolds
\theta}_b)+ g_{22}\bigl(\hat{\bolds\theta}_b^*
\bigr) (\hat\theta_2 - \theta_{0,2}),
\end{equation}
where $\hat{\bolds\theta}_b^*=(\hat\theta_1,\theta_{0,2}^*,\hat\theta_3)'$
satisfies $|\theta_{0,2}^*-\hat\theta_2|\leq|\theta_{0,2}-\hat\theta
_2|$. Therefore, for (\ref{thmtheta2iterat1}) to hold,
it suffices to show that
\begin{eqnarray*}
n^{-(1+\delta)/4}g_2(\hat{\bolds\theta}_b) &\displaystyle\mathop{
\rightarrow}^{d} & N\bigl(0,2^{-1/2}\theta_{0,1}^{-1/2}
\theta_{0,2}^{-3/2}\bigr)\qquad\mbox{if }\delta\in[0,1/3),
\\
g_2(\hat{\bolds\theta}_b) &=& O_p
\bigl(n^\delta\bigr)\qquad \mbox{if }\delta\in[1/3,1),
\\
g_{22}\bigl(\hat{\bolds\theta}_b^*\bigr) &=&
\frac{n^{(1+\delta)/2}} {
2^{3/2}\theta_{0,1}^{1/2}\theta_{0,2}^{3/2}} + o_p\bigl(n^{(1+\delta
)/2}\bigr),
\end{eqnarray*}
where the first two equations follow from (\ref{thmg2order})
with $r_1=(1-\delta)/2$, (\ref{thmg2asymptoticnormality}) and (\ref
{thmtheta1iterate1}), and
the last one is ensured by (\ref{lemmatheta2consistency}), (\ref
{thmg22order}) and
(\ref{thmtheta1iterate1}).
\end{pf*}

\begin{pf*}{Proof of (\protect\ref{thmtheta3iterate1})}
Taking the Taylor expansion of
$g_3(\hat{\bolds\theta})$ at $\hat{\bolds\theta}_c$ yields
%
\begin{eqnarray}\label{thmtheta3taylor}
&& 0 = g_3(\hat{\bolds\theta}) = g_3(\hat{\bolds
\theta}_c)+g_{33}\bigl(\hat{\bolds\theta}_c^*
\bigr) (\hat\theta_3-\theta_{0,3}),
\end{eqnarray}
where $\hat{\bolds\theta}_c^*=(\hat\theta_1,\hat\theta
_2,\theta_{0,3}^*)'$ satisfies
$|\theta_{0,3}^*-\hat\theta_3|\leq|\theta_{0,3}-\hat\theta_3|$.
Therefore, for (\ref{thmtheta3iterate1}) to hold,
it suffices to show that
\begin{eqnarray*}
n^{-\delta/2}g_3(\hat{\bolds\theta}_c)& \displaystyle\mathop{
\rightarrow}^{d}& N\bigl(0,2\theta_{0,3}^{-1}\bigr)\qquad
\mbox{if }\delta\in(0,1/2),
\\
g_3(\hat{\bolds\theta}_c) &=& O_p
\bigl(n^{-(1-3\delta)/2}\bigr)\qquad \mbox{if }\delta\in[1/2,1),
\\
g_{33}\bigl(\hat{\bolds\theta}_c^*\bigr) &=&
\frac{n^\delta}{\theta_{0,3}} + o_p\bigl(n^\delta\bigr),
\end{eqnarray*}
where the first two equations follow from (\ref{thmg3order})
with $r_1=r_2=(1-\delta)/2$, (\ref{thmg3asymptoticnormality}), (\ref
{thmtheta1iterate1}) and~(\ref{thmtheta2iterat1}),
and the last one is ensured by (\ref{thmg33order}).
Thus, (\ref{thmtheta3iterate1}) is\vspace*{-3pt} established.
\end{pf*}

For the second iteration, we show that\vspace*{-3pt}
%
\begin{eqnarray}
n^{1/2} (\hat{\theta}_1-\theta_{0,1} )& \displaystyle\mathop{
\rightarrow}^{d} & N\bigl(0,2\theta_{0,1}^2\bigr)\qquad
\mbox{if }\delta\in[0,1/2),
\nonumber
\\[-8pt]
\label{thmtheta1iterate2}
\\[-8pt]
\nonumber
\hat{\theta}_1-\theta_{0,1} &= & O_p
\bigl(n^{-(1-\delta)}\bigr)\qquad \mbox{if }\delta\in[1/2,1),
\\
n^{(1+\delta)/4} (\hat{\theta}_2-\theta_{0,2} )& \displaystyle\mathop{
\rightarrow}^{d} & N \bigl(0, 2^{5/2}\theta_{0,1}^{1/2}
\theta_{0,2}^{3/2} \bigr)\qquad \mbox{if }\delta\in[0,3/5),
\nonumber
\\[-8pt]
\label{thmtheta2iterat2}
\\[-8pt]
\nonumber
\hat\theta_2 - \theta_{0,2} &=& O_p
\bigl(n^{-(1-\delta)}\bigr)\qquad \mbox{if }\delta\in[3/5,1),
\\
n^{\delta/2} (\hat{\theta}_3-\theta_{0,3} ) &\displaystyle\mathop{
\rightarrow}^{d} & N(0,2\theta_{0,3})\qquad \mbox{if }\delta
\in(0,2/3),
\nonumber
\\[-8pt]
\label{thmtheta3iterate2}
\\[-8pt]
\nonumber
\hat{\theta}_3-\theta_{0,3} &= & O_p
\bigl(n^{-(1-\delta)}\bigr)\qquad \mbox{if }\delta\in[2/3,1).
\end{eqnarray}
By (\ref{thmg1order})
with $r_2=r_3=(1-\delta)/2$, (\ref{thmg1asymptoticnormality}) and (\ref
{thmtheta2iterat1}), we have
\begin{eqnarray*}
n^{-1/2}g_1(\hat{\bolds\theta}_a)& \displaystyle\mathop{
\rightarrow}^{d} & N\bigl(0,2\theta_{0,1}^{-2}
\bigr)\qquad \mbox{if }\delta\in[0,1/2),
\\
g_1(\hat{\bolds\theta}_a) &= & O_p
\bigl(n^\delta\bigr)\qquad \mbox{if }\delta\in[1/2,1).
\end{eqnarray*}
The above two equations,
(\ref{thmg11order}) and (\ref{thmtheta1taylor})
give (\ref{thmtheta1iterate2}). By (\ref{thmg2order})
with $r_1=1-\delta$ and $r_3=(1-\delta)/2$, (\ref
{thmg2asymptoticnormality}), (\ref{thmtheta3iterate1}) and
(\ref{thmtheta1iterate2}), we have
\begin{eqnarray*}
n^{-(1+\delta)/4}g_2(\hat{\bolds\theta}_b)& \displaystyle\mathop{
\rightarrow}^{d} & N\bigl(0,2^{-1/2}\theta_{0,1}^{-1/2}
\theta_{0,2}^{-3/2}\bigr)\qquad \mbox{if }\delta\in[0,3/5),
\\
g_2(\hat{\bolds\theta}_b) &=& O_p
\bigl(n^{-(1-3\delta)/2}\bigr)\qquad \mbox{if }\delta\in[3/5,1).
\end{eqnarray*}
Combining these two equations together with
(\ref{thmg22order}) and (\ref{thmtheta2taylor}) yields (\ref
{thmtheta2iterat2}). By
(\ref{thmg3order})
with $r_1=r_2=1-\delta$, (\ref{thmg3asymptoticnormality}), (\ref
{thmtheta1iterate2}) and (\ref{thmtheta2iterat2}), we
have
\begin{eqnarray*}
n^{-\delta/2}g_3(\hat{\bolds\theta}_c)&  \displaystyle\mathop{
\rightarrow}^{d}& N\bigl(0,2\theta_{0,3}^{-1}\bigr)\qquad
\mbox{if }\delta\in(0,2/3),
\\
g_3(\hat{\bolds\theta}_c) &=& O_p
\bigl(n^{-(1-3\delta)/2}\bigr)\qquad \mbox{if }\delta\in[2/3,1),
\end{eqnarray*}
which, together with (\ref{thmg33order}) and (\ref{thmtheta3taylor}),
lead immediately to
(\ref{thmtheta3iterate2}).

Following the same argument as in the second iteration,
we can recursively show that for each $i=3,4,\ldots$
\begin{eqnarray*}
n^{1/2} (\hat{\theta}_1-\theta_{0,1} )& \displaystyle\mathop{
\rightarrow}^{d} & N\bigl(0,2\theta_{0,1}^2\bigr)\qquad
\mbox{if }\delta\in\bigl[0,(i-1)/i\bigr),
\\
\hat{\theta}_1-\theta_{0,1} &= & O_p
\bigl(n^{-i(1-\delta)/2}\bigr)\qquad \mbox{if }\delta\in\bigl[(i-1)/i,1\bigr),
\\
n^{(1+\delta)/4} (\hat{\theta}_2-\theta_{0,2} )& \displaystyle\mathop{
\rightarrow}^{d} & N \bigl(0, 2^{5/2}\theta_{0,1}^{1/2}
\theta_{0,2}^{3/2} \bigr)\qquad \mbox{if }\delta\in\bigl[0,(2i-1)/(2i+1)\bigr),
\\
\hat\theta_2 - \theta_{0,2} &=& O_p
\bigl(n^{-i(1-\delta)/2}\bigr)\qquad \mbox{if }\delta\in\bigl
[(2i-1)/(2i+1),1\bigr],
\\
n^{\delta/2} (\hat{\theta}_3-\theta_{0,3} )& \displaystyle\mathop{
\rightarrow}^{d} & N(0,2\theta_{0,3})\qquad \mbox{if }\delta\in
\bigl(0,i/(i+1)\bigr),
\\
\hat{\theta}_3-\theta_{0,3} &= & O_p
\bigl(n^{-i(1-\delta)/2}\bigr)\qquad \mbox{if }\delta\in\bigl[i/(i+1),1\bigr).
\end{eqnarray*}
Thus (\ref{thmtheta1consistency})--(\ref
{thmtheta3consistency}) are proved.

\subsection{Proof of Theorem
\protect\ref{thmctlunderincorrectmodel}}\label{sec43}

We divide the proof into three parts corresponding to $\delta\in
[0,1/3)$, $\delta\in[1/3,1/2)$ and $\delta\in[1/2,1)$.

First, we consider $\delta\in[0,1/3)$. We further divide the proof into
six subparts with respect
to $\xi$ in terms of a partition
of $[0,1)$, corresponding to
$\xi\in[0,\delta/2)$, $\xi\in[\delta/2, \delta)$, $\xi\in[\delta
,(1+\delta)/4)$, $\xi\in[(1+\delta)/4,1/2)$,
$\xi\in[1/2,(1+\delta)/2)$ and $\xi\in[(1+\delta)/2,1)$. We shall prove
each of the following six subparts separately:
\begin{longlist}[(a6)]
\item[(a1)] For $\xi\in[(1+\delta)/2,1)$,
\begin{eqnarray*}
&& \hat\theta_1 - \theta_{0,1} = O_p
\bigl(n^{\xi-1}\bigr). 
\end{eqnarray*}
\item[(a2)] For $\xi\in[1/2,(1+\delta)/2)$,
%
\begin{eqnarray}\label{thm2partI-2-1}
\hat\theta_1-\theta_{0,1}& =& O_p
\bigl(n^{\xi-1}\bigr),
\\
\label{thm2partI-2-2}
\hat\theta_2-\theta_{0,2}& =& O_p
\bigl(n^{\xi-(1+\delta)/2}\bigr).
\end{eqnarray}
\item[(a3)] For $\xi\in[(1+\delta)/4,1/2)$,
%
\begin{eqnarray}\label{thm2partI-3-1}
n^{1/2} (\hat\theta_1-\theta_{0,1} )& \displaystyle\mathop{
\rightarrow}^{d}& N\bigl(0,2\theta_{0,1}^2\bigr),
\\
\label{thm2partI-3-2}
\hat\theta_2-\theta_{0,2} &=& O_p
\bigl(n^{\xi-(1+\delta)/2}\bigr).
\end{eqnarray}
\item[(a4)] For $\xi\in[\delta,(1+\delta)/4)$,
%
\begin{eqnarray}\label{thm2partI-4-1}
n^{1/2} (\hat\theta_1-\theta_{0,1} )& \displaystyle\mathop{
\rightarrow}^{d}& N\bigl(0,2\theta_{0,1}^2\bigr),
\\
\label{thm2partI-4-2}
n^{(1+\delta)/4} (\hat\theta_2-\theta_{0,2} )& \displaystyle\mathop{
\rightarrow}^{d}& N\bigl(0,2^{5/2}\theta_{0,1}^{1/2}
\theta_{0,2}^{3/2}\bigr).
\end{eqnarray}
\item[(a5)] For $\xi\in[\delta/2,\delta)$,
%
\begin{eqnarray}\label{thm2partI-5-1}
n^{1/2} (\hat\theta_1-\theta_{0,1} ) & \displaystyle\mathop{
\rightarrow}^{d}& N\bigl(0,2\theta_{0,1}^2\bigr),
\\
\label{thm2partI-5-2}
n^{(1+\delta)/4} (\hat\theta_2-\theta_{0,2} ) &\displaystyle\mathop{
\rightarrow}^{d}& N\bigl(0,2^{5/2}\theta_{0,1}^{1/2}
\theta_{0,2}^{3/2}\bigr),
\\
\label{thm2partI-5-3}
\hat\theta_3-\theta_{0,3}& =& O_p
\bigl(n^{\xi-\delta}\bigr).
\end{eqnarray}
\item[(a6)] For $\xi\in[0,\delta/2)$,
%
\begin{eqnarray}\label{thm2partI-6-1}
n^{1/2} (\hat\theta_1-\theta_{0,1} )& \displaystyle\mathop{
\rightarrow}^{d} & N\bigl(0,2\theta_{0,1}^2\bigr),
\\
\label{thm2partI-6-2}
n^{(1+\delta)/4} (\hat\theta_2-\theta_{0,2} ) &\displaystyle\mathop{
\rightarrow}^{d} & N\bigl(0,2^{5/2}\theta_{0,1}^{1/2}
\theta_{0,2}^{3/2}\bigr),
\end{eqnarray}
and if in addition $\delta\neq0$, then
%
\begin{eqnarray}\label{thm2partI-6-3}
n^{\delta/2} (\hat\theta_3-\theta_{0,3} ) & \displaystyle\mathop{
\rightarrow}^{d} & N(0,2\theta_{0,3}).
\end{eqnarray}
\end{longlist}
\begin{pf*}{Proof of (a1)}
Applying (\ref{thmg1order}) with
$r_1=r_2=r_3=0$ and $\xi\in[(1+\delta)/2,1)$, we have
%
\begin{eqnarray}\label{thm2g1order}
g_1(\hat{\bolds\theta}_a) &=& O_p
\bigl(n^\xi\bigr).
\end{eqnarray}
According to (\ref{lemmatheta1consistency}) and
(\ref{thmg11order}), we have
%
\begin{eqnarray} \label{thm2g11order}
g_{11}(\hat{\bolds\theta}_a) &=& \frac{n}{\theta_{0,1}^2} +
o_p(n).
\end{eqnarray}
The desired conclusion (a1) now follows from plugging
(\ref{thm2g1order}) and (\ref{thm2g11order}) into (\ref{thmtheta1taylor}).
\end{pf*}
\begin{pf*}{Proof of (a2)}
Applying (\ref{thmg1order}) with
$r_1=r_2=r_3=0$ and $\xi\in[1/2,(1+\delta)/2)$, we have
\begin{eqnarray*}
&& g_1(\hat{\bolds\theta}_a)= O_p
\bigl(n^{(1+\delta)/2}\bigr).
\end{eqnarray*}
Combining this with (\ref{thmtheta1taylor}) and (\ref
{thm2g11order}) gives
\begin{eqnarray*}
&& \hat\theta_1-\theta_{0,1} = O_p
\bigl(n^{-(1-\delta)/2}\bigr). 
\end{eqnarray*}
Applying (\ref{thmg2order}) with
$r_1=(1-\delta)/2$, $r_2=r_3=0$ and $\xi\in[1/2,(1+\delta)/2)$, we
obtain
%
\begin{eqnarray}\label{thm2g2order}
g_2(\hat{\bolds\theta}_b) &=& O_p
\bigl(n^\xi\bigr).
\end{eqnarray}
From (\ref{lemmatheta2consistency}) and
(\ref{thmg22order}), we have
%
\begin{eqnarray}\label{thm2g22order}
g_{22}\bigl(\hat{\bolds\theta}_b^*\bigr) &=&
\frac{n^{(1+\delta)/2}}{2^{3/2}\theta_{0,1}^{1/2}\theta_{0,2}^{3/2}} +
o_p\bigl(n^{(1+\delta)/2}\bigr).
\end{eqnarray}
Combining this with (\ref{thmtheta2taylor}) and
(\ref{thm2g2order}) leads to (\ref{thm2partI-2-2}). In addition,
applying (\ref{thmg1order}) with
$r_2=(1+\delta)/2-\xi,r_3=0$ and $\xi\in[1/2,(1+\delta)/2)$, we
have
\begin{eqnarray*}
&& g_1(\hat{\bolds\theta}_a) = O_p
\bigl(n^{1/2}\bigr) + O_p\bigl(n^\xi\bigr) =
O_p\bigl(n^\xi\bigr).
\end{eqnarray*}
This together with (\ref{thmtheta1taylor}) and (\ref{thm2g11order})
gives (\ref{thm2partI-2-1}).
\end{pf*}
\begin{pf*}{Proof of (a3)}
Following the same arguments as
the one used in the proof of
(\ref{thm2partI-2-2}) leads to (\ref{thm2partI-3-2}).
Applying (\ref{thmg1asymptoticnormality}) with
$r_1=(1-\delta)/2$, $r_2=(1+\delta)/2-\xi$, $r_3=0$ and
$\xi\in[(1+\delta)/4,1/2)$, we have
\begin{eqnarray*}
n^{-1/2}g_1(\hat{\bolds\theta}_a)& \displaystyle\mathop{
\rightarrow}^{d}& N\bigl(0,2\theta_{0,1}^{-2}\bigr).
\end{eqnarray*}
This together with (\ref{thmtheta1taylor}) and (\ref{thm2g11order})
gives (\ref{thm2partI-3-1}).
\end{pf*}
\begin{pf*}{Proof of (a4)}
Applying (\ref{thmg2asymptoticnormality})
with $r_1=(1-\delta)/2,r_2=r_3=0$ and
$\xi\in[\delta,(1+\delta)/4)$, we have
\begin{eqnarray*}
&& n^{-(1+\delta)/4}g_2(\hat{\bolds\theta}_b) \displaystyle\mathop{
\rightarrow}^{d} N\bigl(0,2^{-1/2}\theta_{0,1}^{-1/2}
\theta_{0,2}^{-3/2}\bigr).
\end{eqnarray*}
This, (\ref{thmtheta2taylor}) and
(\ref{thm2g22order}) imply (\ref{thm2partI-4-2}). Moreover,
(\ref{thm2partI-4-1}) can be shown by an argument similar to that used
to prove
(\ref{thm2partI-3-1}).
\end{pf*}
\begin{pf*}{Proof of (a5)}
The proofs of
(\ref{thm2partI-5-1}) and (\ref{thm2partI-5-2})
are similar to those of
(\ref{thm2partI-4-1}) and (\ref{thm2partI-4-2}), respectively.
Applying (\ref{thmg3order}) with $r_1=r_2=(1-\delta)/2$,
$r_3=0$ and $\xi\in[(1+\delta)/4,1/2)$, we have
%
\begin{eqnarray}\label{thm2g3order}
g_3(\hat{\bolds\theta}_c)& =& O_p
\bigl(n^\xi\bigr).
\end{eqnarray}
From (\ref{lemmatheta1consistency})--(\ref
{lemmatheta3consistency}) and (\ref{thmg33order}), we obtain
%
\begin{eqnarray}\label{thm2g33order}
&& g_{33}\bigl(\hat{\bolds\theta}_c^*\bigr) =
\frac{n^\delta}{\theta_{0,3}} + o\bigl(n^\delta\bigr).
\end{eqnarray}
Combining this with (\ref{thmtheta3taylor}) and (\ref
{thm2g3order}) leads to (\ref{thm2partI-5-3}).
\end{pf*}
\begin{pf*}{Proof of (a6)}
Equations
(\ref{thm2partI-6-1}) and (\ref{thm2partI-6-2})
can be proved in a way similar to the proofs of
(\ref{thm2partI-4-1}) and (\ref{thm2partI-4-2}).
Applying (\ref{thmg3asymptoticnormality}) with
$r_1=r_2=(1-\delta)/2$, $r_3=0$ and $\xi\in(0,\delta/2)$, we have
\begin{eqnarray*}
n^{-\delta/2}g_3(\hat{\bolds\theta}_c)& \displaystyle\mathop{
\rightarrow}^{d}& N\bigl(0,2{\theta_{0,3}^{-1}}
\bigr).
\end{eqnarray*}
This together with (\ref{thmtheta3taylor}) and (\ref{thm2g33order})
gives (\ref{thm2partI-6-3}).
\end{pf*}

Second, we consider $\delta\in[1/3,1/2)$. Following an argument similar
to that used in the first part, we obtain
\begin{longlist}[(b6)]
\item[(b1)] For $\xi\in[(1+\delta)/2,1)$,
\[
\hat\theta_1 - \theta_{0,1} = O_p
\bigl(n^{\xi-1}\bigr).
\]
\item[(b2)] For $\xi\in[1/2,(1+\delta)/2)$,
\begin{eqnarray*}
\hat\theta_1-\theta_{0,1} &=& O_p
\bigl(n^{\xi-1}\bigr),
\\
\hat\theta_2-\theta_{0,2}& =& O_p
\bigl(n^{\xi-(1+\delta)/2}\bigr).
\end{eqnarray*}
\item[(b3)] For $\xi\in[\delta,1/2)$,
\begin{eqnarray*}
n^{1/2} (\hat\theta_1-\theta_{0,1} )& \displaystyle\mathop{
\rightarrow}^{d}& N\bigl(0,2\theta_{0,1}^2\bigr),
\\
\hat\theta_2-\theta_{0,2} &=& O_p
\bigl(n^{\xi-(1+\delta)/2}\bigr).
\end{eqnarray*}
\item[(b4)] For $\xi\in[(1+\delta)/4,\delta)$,
\begin{eqnarray*}
n^{1/2} (\hat\theta_1-\theta_{0,1} ) &\displaystyle\mathop{
\rightarrow}^{d}& N\bigl(0,2\theta_{0,1}^2\bigr),
\\
\hat\theta_2-\theta_{0,2} &=& O_p
\bigl(n^{\xi-(1+\delta)/2}\bigr),
\\
\hat\theta_3-\theta_{0,3}& =& O_p
\bigl(n^{\xi-\delta}\bigr).
\end{eqnarray*}
\item[(b5)] For $\xi\in[\delta/2,(1+\delta)/4)$,\vspace*{-3pt}
\begin{eqnarray*}
n^{1/2} (\hat\theta_1-\theta_{0,1} ) &\displaystyle\mathop{
\rightarrow}^{d}& N\bigl(0,2\theta_{0,1}^2\bigr),
\\
n^{(1+\delta)/4} (\hat\theta_2-\theta_{0,2} )& \displaystyle\mathop{
\rightarrow}^{d}& N\bigl(0,2^{5/2}\theta_{0,1}^{1/2}
\theta_{0,2}^{3/2}\bigr),
\\
\hat\theta_3-\theta_{0,3} &=& O_p
\bigl(n^{\xi-\delta}\bigr).
\end{eqnarray*}
\item[(b6)] For $\xi\in[0,\delta/2)$,\vspace*{-3pt}
\begin{eqnarray*}
n^{1/2} (\hat\theta_1-\theta_{0,1} ) &\displaystyle\mathop{
\rightarrow}^{d}& N\bigl(0,2\theta_{0,1}^2\bigr),
\\
n^{(1+\delta)/4} (\hat\theta_2-\theta_{0,2} ) &\displaystyle\mathop{
\rightarrow}^{d}& N\bigl(0,2^{5/2}\theta_{0,1}^{1/2}
\theta_{0,2}^{3/2}\bigr),
\\
n^{\delta/2} (\hat\theta_3-\theta_{0,3} )& \displaystyle\mathop{
\rightarrow}^{d}& N(0,2\theta_{0,3}).
\end{eqnarray*}
\end{longlist}

Third, for $\delta\in[1/2,1)$, one can similarly show that
\begin{longlist}[(c6)]
\item[(c1)] For $\xi\in[(1+\delta)/2,1)$,
\[
\hat\theta_1 - \theta_{0,1} = O_p
\bigl(n^{\xi-1}\bigr).
\]
\item[(c2)] For $\xi\in[\delta,(1+\delta)/2)$,
\begin{eqnarray*}
\hat\theta_1-\theta_{0,1} &=& O_p
\bigl(n^{\xi-1}\bigr),
\\
\hat\theta_2-\theta_{0,2}& =& O_p
\bigl(n^{\xi-(1+\delta)/2}\bigr).
\end{eqnarray*}
\item[(c3)] For $\xi\in[1/2,\delta)$,
\begin{eqnarray*}
\hat\theta_1-\theta_{0,1} &=& O_p
\bigl(n^{\xi-1}\bigr),
\\
\hat\theta_2-\theta_{0,2} &=& O_p
\bigl(n^{\xi-(1+\delta)/2}\bigr),
\\
\hat\theta_3-\theta_{0,3} &=& O_p
\bigl(n^{\xi-\delta}\bigr).
\end{eqnarray*}
\item[(c4)] For $\xi\in[(1+\delta)/4,1/2)$,
\begin{eqnarray*}
n^{1/2} (\hat\theta_1-\theta_{0,1} ) &\displaystyle\mathop{
\rightarrow}^{d}& N\bigl(0,2\theta_{0,1}^2\bigr),
\\
\hat\theta_2-\theta_{0,2} &=& O_p
\bigl(n^{\xi-(1+\delta)/2}\bigr),
\\
\hat\theta_3-\theta_{0,3}& =& O_p
\bigl(n^{\delta/2}\bigr).
\end{eqnarray*}
\item[(c5)] For $\xi\in[\delta/2,(1+\delta)/4)$,
\begin{eqnarray*}
n^{1/2} (\hat\theta_1-\theta_{0,1} ) &\displaystyle\mathop{
\rightarrow}^{d}& N\bigl(0,2\theta_{0,1}^2\bigr),
\\
n^{(1+\delta)/4} (\hat\theta_2-\theta_{0,2} )& \displaystyle\mathop{
\rightarrow}^{d}& N\bigl(0,2^{5/2}\theta_{0,1}^{1/2}
\theta_{0,2}^{3/2}\bigr),
\\
\hat\theta_3-\theta_{0,3} &=& O_p
\bigl(n^{\xi-\delta}\bigr).
\end{eqnarray*}
\item[(c6)] For $\xi\in[0,\delta/2)$,
\begin{eqnarray*}
n^{1/2} (\hat\theta_1-\theta_{0,1} ) &\displaystyle\mathop{
\rightarrow}^{d}& N\bigl(0,2\theta_{0,1}^2\bigr),
\\
n^{(1+\delta)/4} (\hat\theta_2-\theta_{0,2} )& \displaystyle\mathop{
\rightarrow}^{d}& N\bigl(0,2^{5/2}\theta_{0,1}^{1/2}
\theta_{0,2}^{3/2}\bigr),
\\
n^{\delta/2} (\hat\theta_3-\theta_{0,3} )& \displaystyle\mathop{
\rightarrow}^{d}& N(0,2\theta_{0,3} ).
\end{eqnarray*}
\end{longlist}
Thus the proof of the theorem is complete.

\subsection{Proofs of Corollaries \protect\ref{coropolymle0} and
\protect\ref{coropolymle}}\label{sec44}
To prove Corollaries \ref{coropolymle0} and \ref{coropolymle}, the
following lemma, which provides the order of magnitude of $R(\Theta)$
defined in
(\ref{eqrisk}), is needed.
%
\begin{lem}
\label{lemmapolyequations}
Under the setup of Lemma~\ref{lem02Ttraces}, let
$\mathbf{x}=n^{-1}(1,2,\dots,n)'$ and $\mathbf{1}=(1,\ldots,1)'$.
Then for any $\delta\in[0,1)$, the following equations hold uniformly
in $\Theta$:
%
\begin{eqnarray}\label{polyresult1}
\mathbf{1}'\bolds\Sigma^{-1}(\bolds\theta)\mathbf{1} &=&
\frac{\theta_3^2}{2\theta_2}n^\delta+ o\bigl(n^\delta\bigr) +O(1),
\\
\label{polyresult2}
\mathbf{x}'\bolds\Sigma^{-1}(\bolds\theta)\mathbf{1} &=&
\frac{\theta_3^2}{4\theta_2}n^\delta+ o\bigl(n^\delta\bigr) + O(1),
\\
\label{polyresult3}
\mathbf{x}'\bolds\Sigma^{-1}(\bolds\theta)\mathbf{x}& =&
\frac{\theta_3^2}{6\theta_2}n^\delta+ o\bigl(n^\delta\bigr) + O(1).
\end{eqnarray}
\end{lem}

We first prove Corollary~\ref{coropolymle0}. Note that
%
\begin{eqnarray}
&&\bolds\mu_0'\bolds\Sigma^{-1}  (\bolds
\theta) \bigl(\mathbf{I}-\mathbf{M}(\bolds\theta)\bigr)\bolds\mu_0
\nonumber\\
&&\quad= \bolds\mu_0'\bolds\Sigma^{-1}(\bolds
\theta)\bolds\mu_0 - \bolds\mu_0'\bolds
\Sigma^{-1}(\bolds\theta)\mathbf{M}(\bolds\theta)\bolds\mu_0
\nonumber
\\[-8pt]
\label{cor1misselectorder}
\\[-8pt]
\nonumber
&&\quad= \beta_{0,1}^2\mathbf{x}'\bolds
\Sigma^{-1}(\bolds\theta)\mathbf{x} -\beta_{0,1}^2
\mathbf{x}'\bolds\Sigma^{-1}(\bolds\theta)\mathbf{1}\bigl(
\mathbf{1}'\bolds\Sigma^{-1}(\bolds\theta)\mathbf{1}
\bigr)^{-1} \mathbf{1}'\bolds\Sigma^{-1}(\bolds
\theta)\mathbf{x}
\\
&&\quad= \frac{\beta_{0,1}^2\theta_3^2}{24\theta_2}n^\delta+ o\bigl(n^\delta
\bigr),\nonumber
\end{eqnarray}
uniformly in $\Theta$, where $\mathbf{x}=n^{-1}(1,\dots,n)'$
and the last equality is obtained from (\ref{polyresult1})--(\ref
{polyresult3}). Therefore, (\ref{cor1xi=theta}) holds.
With the help of (\ref{cor1xi=theta}),
(\ref{cor1theta1consistency}) and
(\ref{cor1theta2consistency}) follow directly from Theorem~\ref{thmctlunderincorrectmodel}.

Second, we prove Corollary~\ref{coropolymle}. By
(\ref{loglikedecompose}),
(\ref{lemmamis2}) and (\ref{cor1misselectorder}), we have
\begin{eqnarray*}
-2\ell(\bolds\theta) &=& n\log(2\pi)+\log\det\bigl(\bolds\Sigma(\bolds
\theta)
\bigr) +\operatorname{tr}\bigl(\bolds\Sigma^{-1}(\bolds\theta)\bolds\Sigma(
\bolds\theta_0)\bigr)
\\
&&{}+ \frac{\beta_{0,1}^2\theta_3^2}{24\theta_2}n^\delta+ h(\bolds\theta) +o_p
\bigl(n^\delta\bigr)+ O_p(1),
\end{eqnarray*}
uniformly in $\Theta$, noting that $h(\bolds\theta)= (\bolds
\eta+\bolds\epsilon)'\bolds\Sigma^{-1}(\bolds\theta)(\bolds\eta+\bolds
\epsilon)
-\operatorname{tr}(\bolds\Sigma^{-1}(\bolds\theta)\bolds\Sigma(\bolds\theta
_0))$. Therefore,
by (\ref{logdetSigma}) and (\ref{Sigma1Sigmainverse3}),
%
\begin{eqnarray}
-2\ell(\bolds\theta) &=& n\log(2\pi) - \frac{1-\delta}{2}\log{n} + \biggl
(\log
\theta_1 + \frac{\theta_{0,1}}{\theta_1} \biggr)n
\nonumber\\
&&{}+ \biggl(\frac{2\theta_2}{\theta_1} \biggr)^{1/2} \biggl(1-\frac
{\theta_{0,1}}{2\theta_1}
+\frac{\theta_{0,2}}{2\theta_2} \biggr)n^{(1+\delta)/2}
\nonumber\\
&&{}- \biggl(\frac{\theta_2}{\theta_1} +\theta_3-\frac{\theta_{0,2}(\theta
_3^2 - \theta_{0,3}^2)}{2\theta_2\theta_{0,3}} -
\frac{\beta_{0,1}^2\theta_3^2}{24\theta_2} \biggr)n^\delta
\nonumber\\
&&{}+ h(\bolds\theta)+o_p\bigl(n^\delta\bigr)+
O_p(1)
\nonumber
\\[-8pt]
\label{loglikedecomposepoly}
\\[-8pt]
\nonumber
&=& n\log(2\pi) - \frac{1-\delta}{2}\log{n} + \biggl(\log\theta_1 +
\frac{\theta_{0,1}}{\theta_1} \biggr)n
\\
&& {}+ \biggl(\frac{2\theta_2}{\theta_1} \biggr)^{1/2} \biggl(1-\frac
{\theta_{0,1}}{2\theta_1}
+\frac{\theta_{0,2}}{2\theta_2} \biggr)n^{(1+\delta)/2}\nonumber
\\
&&{}- \biggl\{\frac{\theta_2}{\theta_1}+\theta_3 \biggl(1-\frac{\theta
_{0,2}}{\theta_2}
\biggr) +\frac{\theta_{0,2}\theta_{0,3}+\theta_{0,2}\theta
_{0,3}^*}{2\theta_2}\nonumber
\\
&&{}-\frac{\theta_{0,2}}{2\theta_2\theta_{0,3}^*} \bigl(\theta_3-\theta
_{0,3}^*
\bigr)^2 \biggr\} n^\delta+h(\bolds\theta)+o_p
\bigl(n^\delta\bigr)+O_p(1),\nonumber
\end{eqnarray}
uniformly in $\Theta$, where $ \theta_{0,3}^*=
\frac{12\theta_{0,2}}{
12\theta_{0,2}+\beta_{0,1}^2\theta_{0,3}}\theta_{0,3}$.
It follows from
(\ref{loglikedecomposepoly})
and the same argument as in the proof of (\ref{theta3consistent})
that for any $\varepsilon_3>0$, there exist $\varepsilon_1,\varepsilon
_2>0$ such that
\begin{eqnarray*}
&& P \Bigl(\inf_{\bolds\theta\in\Theta_2(\varepsilon_1,\varepsilon
_2,\varepsilon_3)} \bigl\{-2\ell(\bolds\theta)+2\ell\bigl(
\bigl(\theta_1,\theta_2,\theta_{0,3}^*
\bigr)'\bigr)\bigr\}>0 \Bigr)\rightarrow1,
\end{eqnarray*}
as $n\rightarrow\infty$, where
$\Theta_3(\varepsilon_1,\varepsilon_2,\varepsilon_3)
= \{\bolds\theta\in\Theta:|\theta_1-\theta_{0,1}|\leq\varepsilon
_1,|\theta_2-\theta_{0,2}|\leq
\varepsilon_2,|\theta_3-\theta_{0,3}^*|>\varepsilon_3 \}$.
Thus
(\ref{cor1theta3inconsistent}) is established, and hence the proof is complete.

\subsection{Proofs of Corollaries \protect\ref{coroexp0} and \protect
\ref{coroexp}}\label{sec45}
We first prove (\ref{cor1xi=1+theta2}).
Let $\mathbf{x}=(x(s_1),\dots,x(s_n))'$.
By an argument similar to that
used to prove \eqref{lemmamis2},
it can be shown that\vspace*{-6pt}
\begin{eqnarray*}
&&\sup_{\bolds\theta\in\Theta}n^{-\delta/2}\mathbf{x}'\bolds
\Sigma^{-1}(\bolds\theta)\mathbf{1} = O_p(1).
\end{eqnarray*}
This, together with
(\ref{polyresult3}) and
(\ref{Sigma1Sigmainverse}), gives
\begin{eqnarray*}
&&\bolds\mu_0'\bolds\Sigma^{-1}  (\bolds
\theta) \bigl(\mathbf{I}-\mathbf{M}(\bolds\theta)\bigr)\bolds\mu_0
\nonumber\\
&&\quad= \bolds\mu_0'\bolds\Sigma^{-1}(\bolds
\theta)\bolds\mu_0 - \bolds\mu_0'\bolds
\Sigma^{-1}(\bolds\theta)\mathbf{M}(\bolds\theta)\bolds\mu_0
\nonumber\\
&&\quad= \beta_{0,1}^2\mathbf{x}'\bolds
\Sigma^{-1}(\bolds\theta)\mathbf{x} -\beta_{0,1}^2
\mathbf{x}'\bolds\Sigma^{-1}(\bolds\theta)\mathbf{1}\bigl(
\mathbf{1}'\bolds\Sigma^{-1}(\bolds\theta)\mathbf{1}
\bigr)^{-1} \mathbf{1}'\bolds\Sigma^{-1}(\bolds
\theta)\mathbf{x}
\\
&&\quad= \beta_{0,1}^2\operatorname{tr}\bigl(\bolds
\Sigma^{-1}(\bolds\theta)\bolds\Sigma_\eta(0,
\theta_{1,2},\theta_{1,3})'\bigr) +
h_x(\bolds\theta)+ O_p(1)
\nonumber\\
&&\quad= \frac{\beta_{0,1}^2\theta_{1,2}}{(2\theta_1\theta
_2)^{1/2}}n^{(1+\delta)/2} +\frac{\beta_{0,1}^2\theta_{1,2}(\theta
_3^2-\theta_{1,3}^2)}{2\theta_2\theta_{1,3}}n^\delta
\\
&&\qquad{}+ h_x(\bolds\theta)+o\bigl(n^\delta\bigr)+
O_p(1),\nonumber 
\end{eqnarray*}

\noindent uniformly in $\Theta$, where $
h_x(\bolds\theta)=\beta_{0,1}^2 (\mathbf{x}'\bolds\Sigma^{-1}(\bolds
\theta)\mathbf{x}
-\operatorname{tr}(\bolds\Sigma^{-1}(\bolds\theta)\bolds\Sigma_\eta(0,\theta
_{1,2},\theta_{1,3})') )$.
In addition, an argument similar to that used to prove
\eqref{lemmaconvergea.s.5}
yields
\begin{eqnarray*}
&&\sup_{\bolds\theta\in\Theta}h_x(\bolds\theta)= o_p
\bigl(n^{(1+\delta)/2}\bigr).
\end{eqnarray*}
Hence (\ref{cor1xi=1+theta2}) follows.
In view of (\ref{cor1xi=1+theta2}) and Theorem~\ref{thmctlunderincorrectmodel},
we obtain
(\ref{cor3theta1}). Thus, the proof of Corollary~\ref{coroexp0} is complete.

To prove (\ref{cor3theta2inconsistent}), note first that by
the same line of reasoning as in
(\ref{loglikedecomposepoly}), one gets
%
\begin{eqnarray}
-2\ell(\bolds\theta) &=& n\log(2\pi)-\frac{1-\delta}{2}\log n+ \biggl
(\log
\theta_1+\frac{\theta_{0,1}}{\theta_1} \biggr)n
\nonumber\\
& &{}+ \biggl(\frac{2\theta_2}{\theta_1} \biggr)^{1/2} \biggl(1-
\frac{\theta_{0,1}}{2\theta_1}+\frac{\theta_{0,2}^*}{2\theta_2} \biggr
)n^{(1+\delta)/2}
\nonumber\\
\label{loglikedecomposeexp}
&&{}- \biggl\{\frac{\theta_2}{\theta_1}+\theta_3 \biggl(1-\frac{\theta
_{0,2}^*}{\theta_2}
\biggr) +\frac{\theta_{0,2}\theta_{0,3}+\beta_{0,1}^2\theta_{1,2}\theta_{1,3}
+\theta_{0,2}^*\theta_{0,3}^*}{2\theta_2}
\\
&&{}-\frac{\theta_{0,2}^*}{2\theta_2\theta_{0,3}^*} \bigl(
\theta_3 -
\theta_{0,3}^* \bigr)^2 \biggr\}n^\delta
\nonumber\\
&&{}+h_x(\bolds\theta)+h(\bolds\theta)+ o_p
\bigl(n^\delta\bigr) + O_p(1),\nonumber
\end{eqnarray}
uniformly\vspace*{1pt} in $\Theta$, where $ \theta_{0,2}^*=\theta
_{0,2}+\beta_{0,1}^2\theta_{1,2}$ and
$\theta_{0,3}^*= \frac{\theta_{0,2}+\beta_{0,1}^2\theta_{1,2}}{
\beta_{0,1}^2\theta_{1,2}\theta_{1,3}^{-1}+\theta_{0,3}\theta_{0,3}^{-1}}$.
Moreover,
using arguments similar to those used in the proofs of
(\ref{theta2consistent}) and (\ref{theta3consistent}), respectively,
one can show that
for any $\varepsilon_2>0$, there exists an $\varepsilon_1>0$ such that
%
\begin{eqnarray}\label{Ing0507a}
&& \lim_{n \to\infty}P \Bigl(\inf_{\bolds\theta\in\Theta_2(\varepsilon
_1,\varepsilon_2)} \bigl\{-2\ell(
\bolds\theta)+2\ell\bigl(\bigl(\theta_1,\theta_{0,2}^*,
\theta_{0,3}^*\bigr)'\bigr)\bigr\}>0 \Bigr)=1,
\end{eqnarray}
and
for any $\varepsilon_3>0$, there exist
$\varepsilon_1,\varepsilon_2>0$ such that
%
\begin{eqnarray}\label{Ing0507b}
&& \lim_{n \to\infty} P \Bigl(\inf_{\bolds\theta\in\Theta_3(\varepsilon
_1,\varepsilon_2,\varepsilon_3)} \bigl\{-2
\ell(\bolds\theta)+2\ell\bigl(\bigl(\theta_1,\theta_2,
\theta_{0,3}^*\bigr)'\bigr)\bigr\}>0 \Bigr)=1,
\end{eqnarray}
where\vspace*{1pt} $\Theta_2(\varepsilon_1,\varepsilon_2)
= \{\bolds\theta\in\Theta:|\theta_1-\theta_{0,1}|\leq\varepsilon
_1,|\theta_2-\theta_{0,2}^*|>\varepsilon_2 \}$
and
$\Theta_3(\varepsilon_1,\varepsilon_2,\varepsilon_3)
= \{\bolds\theta\in\Theta:|\theta_1-\theta_{0,1}|\leq\varepsilon
_1,|\theta_2-\theta_{0,2}^*|\leq
\varepsilon_2,|\theta_3-\theta_{0,3}^*|>\varepsilon_3 \}$.
Combining\vspace*{1pt}
(\ref{loglikedecomposeexp})--(\ref{Ing0507b}) yields
(\ref{cor3theta2inconsistent})
and~(\ref{cor3theta3inconsistent}).
This completes the proof of Corollary~\ref{coroexp}.

\section*{Acknowledgements}
The authors would like to thank the Associate Editor and an anonymous reviewer
for their insightful and constructive comments, which greatly improve
the presentation
of this paper.
The research of Chih-Hao Chang and Hsin-Cheng Huang was supported by
Ministry of Science and Technology of Taiwan under grants
MOST 103-2118-M-390-005-MY2 and MOST 100-2628-M-001-004-MY3, respectively.
The research of Ching-Kang Ing was supported by
Academia Sinica Investigator Award.

\begin{supplement}
\stitle{Supplement to ``Mixed domain asymptotics for a stochastic process model with
time trend and measurement error''}
\slink[doi]{10.3150/15-BEJ740SUPP} 
\sdatatype{.pdf}
\sfilename{BEJ740\_supp.pdf}
\sdescription{The supplementary material contains the proofs
of lemmas in Section~\ref{sectionproofoftheoremsandcorollaries}, following some technical lemmas
needed in the proofs.}
\end{supplement}

%




\printhistory
\end{document}